%


\magnification = 1200
\raggedbottom

\input amssym.def
\input amssym.tex

\font\bigtenrm=cmr10 scaled\magstep1

\font\eightrm=cmr8

\def\AP{1}  
\def\BHMWZ{2}  
\def\BMW{3}  
\def\BRW{4} 
\def\Ber{5}  
\def\Bera{6}  
\def\BerSid{7}  
\def\DK{8}     
\def\Hal{9}   
\def\Jech{10}  
\def\Kuna{11}  
\def\MW{12}  
\def\RuRu{13}  
\def\Sid{14}    
\def\Sik{15}     

\def\RR{{\Bbb R}}
\def\PP{{\Bbb P}}
\def\QQ{{\Bbb Q}}
\def\NN{{\Bbb N}}
\def\LL{{\Bbb L}}
\def\CC{{\Bbb C}}
\def\BB{{\Bbb B}}
\def\HH{{\Bbb H}}

\def\NNN{{\cal N}}

\def\FF{{\cal F}}

\def\UU{{\cal U}}
\def\WW{{\cal W}}

\def\cc{{\frak c}}

\def\cclw{\cc^{<\omega}}
\def\Dlw{\Delta^{<\omega}}
\def\tlw{2^{<\omega}}

\def\cl{\overline}
\def\imp{\,\Rightarrow\,}
\def\iff{\,\Leftrightarrow\,}  

\def\rest{\upharpoonright}  

\def\eop {{\vrule height4pt width4pt depth0pt}}  

\centerline{\bigtenrm Locally Constant Functions}

\centerline {Joan Hart${}^1$ and Kenneth Kunen\footnote{${}^1$}{\sevenrm Authors
supported by NSF Grant DMS-9100665.}}

\centerline {\sevenrm University of Wisconsin}

\centerline {\sevenrm Madison, WI  53706,\ \ U.S.A.}

\centerline {\sevenrm jhart@math.wisc.edu and kunen@cs.wisc.edu}

\centerline{\eightrm May 5, 1995}

\bigskip
\bigskip

\centerline{\bf ABSTRACT}
{\narrower\medskip\noindent
Let $X$ be a compact Hausdorff space and $M$ a metric space.
$E_0(X,M)$ is the set of 
$f \in C(X,M)$ such that there is a dense set of points $x\in X$
with $f$ constant on some neighborhood of $x$.
We describe some general classes of $X$ for which
$E_0(X,M)$ is all of $C(X,M)$.
These include $\beta\NN \backslash \NN$, any nowhere separable LOTS,
and any $X$ such that forcing with the open subsets of $X$ does not
add reals.
In the case that $M$ is a Banach space, we discuss the properties
of $E_0(X,M)$ as a normed linear space.
We also build three first countable Eberlein compact spaces, $F,G,H$,
with various $E_0$ properties.   For all metric $M$,
$E_0(F,M)$ contains only the constant functions, and $E_0(G,M) = C(G,M)$.
If $M$ is the Hilbert cube or any infinite dimensional Banach space,
$E_0(H,M) \ne C(H,M)$, but $E_0(H,M) = C(H,M)$ whenever $M \subseteq \RR^n$
for some finite $n$.

\bigskip
}

\bigskip

{\bf \S0. Introduction.}  If $X$ is a compact Hausdorff space 
and $M$ is a metric space, let $C(X,M)$ be the space of all continuous
functions from $X$ into $M$.  $C(X,M)$ is a metric space under
the sup norm. $C(X)$ denotes $C(X,\RR)$, which is a (real) Banach algebra.
Following [\Ber, \Bera, \BerSid, \RuRu, \Sid], if
$f \in C(X,M)$, let $\Omega_f$ be the union of
all open $U \subseteq X$ such that $f$ is constant on $U$.
Then, $E_0(X,M)$ is the set of all $f \in C(X,M)$ such
that $\Omega_f$ is dense in $X$; these functions are called ``locally
constant on a dense set''.
$E_0(X)$ denotes $E_0(X,\RR)$.

Clearly, $E_0(X)$ is a subalgebra
of $C(X)$ and contains all the constant functions.
As Bernard and Sidney point out [\Bera, \BerSid, \Sid],
if $X$ is compact metric with no isolated points, then
$E_0(X)$ is a proper dense subspace of $C(X)$.
In this paper, we study the two extreme situations:  where $E_0(X)$ contains
only the constant functions, and where $E_0(X) = C(X)$.
In \S5, we give some justification for studying these two extremes.

A standard example of elementary analysis is a monotonic
$f \in C([0,1])$ which does all its growing on a Cantor set;
then $f$ is a nonconstant function in 
$E_0([0,1])$.  More generally, for ``many'' $X$,
$E_0(X)$ separates points in $X$, and hence (by the
Stone-Weierstrass Theorem), is dense in $C(X)$.  Specifically,

\medskip

{\bf 0.1. Theorem.}  If $X$ is compact Hausdorff and
$E_0(X)$ is {\it not\/} dense in
$C(X)$, then 

\item{a.} $X$ has a family of $2^{\aleph_0}$ disjoint nonempty open subsets.
\item{b.} $X$ is not locally connected.
\item{c.} $X$ is not zero-dimensional.

\medskip

Part (c) of the Theorem is obvious.  Parts (a) and (b) are due
to M. E. Rudin and W. Rudin [\RuRu], and generalize earlier
results of Bernard and Sidney that if $X$ is compact and
second countable, then $E_0(X)$ is  dense in $C(X)$.

However, {\it first\/} countable is not enough.  In \S2, we produce
a first countable compact $X$ such that $E_0(X)$ contains only
the constant functions.  A {\it non\/} first countable example
was constructed in [\RuRu].  Our result is patterned after [\RuRu].
Roughly, we replace the family of Cantor sets used in their construction
by a {\it disjoint\/} family.  This adds some complexity to the
construction.  However, we also simplify the
geometry of the construction by building the space inside a Hilbert space.
Our space will be compact in the weak topology, and hence
a uniform Eberlein compact (that is, a weakly compact
subspace of a Hilbert space).   One may use the approach of \S2
to simplify the construction of [\RuRu] and to demonstrate that their
space is also a uniform Eberlein compact.

In \S3, we look at the other extreme; there are many familiar compact
$X$ for which $E_0(X)$ is all of $C(X)$, such as
$\beta \NN \backslash \NN$ and a 
Suslin line. For some classes of spaces, such as
compact ordered spaces and compact extremally disconnected spaces,
we present simple necessary and sufficient conditions
for $E_0(X) = C(X)$.

In \S\S3,4, we also consider $E_0(X,M)$ for
other metric spaces $M$.  It is easy to see that
$E_0(X,\CC) = C(X,\CC)$ iff $E_0(X,\RR) = C(X,\RR)$,
and 
$E_0(X,\CC)$ is dense in $C(X,\CC)$ iff $E_0(X,\RR)$ is dense in $C(X,\RR)$,
but the situation for general $M$ is a bit more complex.
In particular, in \S4, we produce a
uniform Eberlein compact $X$ such that 
$E_0(X,\RR) = C(X,\RR)$ but
$E_0(X,Q) \ne C(X,Q)$, where $Q$ is the Hilbert cube.
In \S5, we let $M$ be a Banach space, and consider the properties
of $E_0(X,M)$ as a normed linear space.

Also in \S5, we show that $E_0(X)$ is a proper dense subspace of
$C(X)$ whenever $X$ is a nontrivial infinite product.

In \S1, we prove some preliminary results on Cantor sets used 
in our construction in \S2.

Independently of Bernard and Sidney,
\ \ Bella, Hager, Martinez, Woodward, and Zhou [\BHMWZ, \BMW, \MW]
defined the space $E_0(X)$ (they called it $dc(X)$), and
showed (in the spirit of Theorem 0.1) that $E_0(X)$ is dense
in $C(X)$ in many cases.  We comment further on their work at the end of \S3.

\bigskip

{\bf \S1. Cantor Sets.} 
By a {\it closed interval\/} we mean any compact space homeomorphic
to $[0,1] \subseteq \RR$.  By a {\it Cantor set\/} we mean any space
homeomorphic to the usual Cantor set in $\RR$; equivalently, 
homeomorphic to $2^\omega$, where $2 = \{0,1\}$ has the
discrete topology.  The following lemma was used also in [\RuRu].

\medskip

{\bf 1.1. Lemma.}  If $J$ is a closed interval, $f \in C(J)$, and $f$
is not constant, then there is a Cantor set $H \subset J$
such that $f$ is $1 - 1$ on $H$.

\medskip

In our construction, we need a uniform version of this.
If $H$ is a subset of a product $X \times J$, we
use $H_x$ to denote $\{y\in J : (x,y) \in H\}$.

\medskip

{\bf 1.2. Lemma.}  Suppose $J$ is a closed interval and $X$ is
a compact zero-dimensional Hausdorff space, and suppose
$f \in C(X \times J)$ is such that for every $x \in X$,
$f \rest (\{x\} \times J)$ is not constant.  Then there is
a set $H \subset X \times J$ such that:

\vbox{
\item{(1)} $H_x$ is a Cantor set for every $x \in X$.
\item{(2)} $f$ is $1 - 1$ on $\{x\} \times H_x$ for every $x \in X$.
\item{(3)} There is a continuous $\varphi : H \to 2^\omega$ such
that the map $(x,y) \mapsto (x, \varphi(x,y))$ is a homeomorphism
from $H$ onto $X \times 2^\omega$.
}

\medskip

Some remarks:  Lemma 1.1 is the special case of Lemma 1.2 
where $X$ is a singleton.  If we deleted (3), then 1.2 would be immediate
from 1.1, using the Axiom of Choice, without any assumption on $X$.
But (3) says that we can choose the Cantor sets continuously.
As stated, the Theorem requires $X$ to be zero-dimensional.
For example, suppose $X = J = [0,1]$, and we take $f$
to be constant on the strip $\{(x,y) : |x - y| < {1 \over 3}\}$.
Then $H$ must be disjoint from the strip, which is easily
seen to contradict (3).  Of course, (1) follows from (3).

Lemma 1.1 may be proved by a binary tree argument, and we prove
Lemma 1.2 by showing how to build this tree ``uniformly'' for all $x \in X$.
A simpler proof of Lemma 1.1 in [\RuRu] takes
advantage of the ordering on $\RR$, but this proof does not easily generalize
to a proof of Lemma 1.2.  Moreover, the tree argument 
extends to non-ordered spaces.
For example, in Lemma 1.2, $J$ could be any compact metric space
which is connected and locally connected, and $f$ could be any 
map into a Hausdorff space.

The following general tree notation will be used here and in \S\S2-4.
If $\Delta$ is some index set, then $\Dlw$ denotes the  tree of all
finite sequences from $\Delta$; this is the complete
$\Delta$-ary tree of height $\omega$.
For $s \in \Dlw$, 
let $lh(s) \in \omega$ be its length.  We use $()$ to denote the 
empty sequence.
If $i \le lh(s)$, let $s\rest i$ be the sequence of length $i$ consisting
of the first $i$ elements of $s$;
$t \subseteq s$ iff $t = s \rest i$ for some $i \le lh(s)$.
Let $t\alpha$ denote the sequence of length $lh(t) + 1$
obtained by appending $\alpha$ to $t$.
Note that $\Dlw$, ordered by $\subseteq$,
is a tree with root $()$, and the nodes immediately above $s$ are
the $s\alpha$ for $\alpha\in\Delta$.
We say $s,t \in \Dlw$ are 
{\it compatible} iff $s \subseteq t$ or $t \subseteq s$.
We let $s \perp t$ abbreviate the statement that
$s,t$ are {\it incompatible\/} (not compatible).

A {\it path\/} in $\Dlw$ is a chain, $P$, such that
$s \alpha  \in P$ implies $s \in P$ for all $s$ and $\alpha$.
A path may be empty or finite or countably infinite.
The infinite paths are all of the form $\{\psi \rest n : n \in \omega\}$
where $\psi: \omega \to \Delta$.
In particular, for binary trees, $\Delta = 2 = \{0,1\}$,
and the infinite paths through the Cantor tree, $\tlw$, are
associated with the points in the Cantor set, $2^\omega$.

To prove 1.2, fix a metric on $J$. For $E \subseteq J$, let
$diam(E)$ be the diameter of $E$ with respect to this metric.
Call a subset of $X \times J$ {\it simple\/} iff it is of the form
$\bigcup_{i < k}Q_i\times I_i$, where  $k$ is finite,
the $Q_i$ for $i < k$ form a disjoint family of clopen sets whose
union is $X$, and each $I_i$ is a closed interval.  We prove 1.2
by iterating the following splitting lemma.

\medskip

{\bf 1.3. Lemma.} Let $J,X,f$ be as in 1.2 and let $\epsilon > 0$.
Then there are simple $A_0,A_1 \subset X\times J$ such that
the following hold:

\item{(a)} $A_0 \cap A_1 = \emptyset$.
\item{(b)} For each $x \in X$,
$f(\{x\} \times (A_0)_x) \cap f(\{x\} \times (A_1)_x) = \emptyset$.
\item{(c)} For each $x \in X$ and $\mu = 0,1$,
$diam((A_\mu)_x) \le \epsilon$.
\item{(d)} For each $x \in X$ and $\mu = 0,1$,
$f \rest (\{x\} \times (A_\mu)_x)$ is not constant.

{\bf Proof.}  For each $z \in X$, $f \rest (\{z\}  \times J)$ is
a nonconstant map from an interval into an interval, so we may
choose disjoint closed intervals $I_0(z), I_1(z) \subset J$
such that $f(\{z\} \times I_0(z)) \cap f(\{z\} \times I_1(z)) = \emptyset$,
$diam(I_\mu(z))\le \epsilon$, and
$f \rest (\{z\}  \times I_\mu(z))$ is not constant ($\mu = 0,1$).
By continuity, there is a neighborhood $U_z$ of $z$ such that
for all $x \in U_z$, 
$f(\{x\} \times I_0(z)) \cap f(\{x\} \times I_1(z)) = \emptyset$ and
$f \rest (\{x\}  \times I_\mu(z))$ is not constant.
Since $X$ is compact and zero-dimensional, there are a finite $k$,
points $z_i \in X$ ($i < k$), and clopen $Q_i \subseteq U_{z_i}$ such
that the $Q_i$ form a partition of $X$.  Then let 
$A_\mu = \bigcup_{i < k}Q_i\times I_\mu(z_i)$.  \eop

\medskip

\vbox{  
{\bf Proof of 1.2.}  For $s \in \tlw$, choose simple $A_s \subseteq X \times J$
such that 

\item{(a)} For each $s\in\tlw$, $A_{s0} \cap A_{s1} = \emptyset$.
\item{(b)} For each $x \in X$ and $s\in\tlw$ ,
$f(\{x\} \times (A_{s0})_x) \cap f(\{x\} \times (A_{s1})_x) = \emptyset$.
\item{(c)} For each $x \in X$ and $t \in \tlw$,
$diam((A_t)_x) \le 1/lh(t)$.
\item{(d)} For each $x \in X$ and  $t \in \tlw$,
$f \rest (\{x\} \times (A_t)_x)$ is not constant.
}

\noindent
We may take $A_{()} = X \times J$; then, for $t = ()$, (c) is
vacuous and (d) follows from the hypothesis of 1.2.  Given $A_s$,
we obtain $A_{s0}$ and $A_{s1}$ by applying 1.3 to each box making
up $A_s$.  Let $H = \bigcap_{n \in \omega} \bigcup\{A_s: lh(s) = n\}$.
Let $\varphi(x,y)$ be the (unique) $\psi \in 2^\omega$ such that
$(x,y) \in A_{\psi\rest n}$ for all $n \in \omega$.  \eop

\bigskip

{\bf \S2. Making $E_0(X)$ Small.}  We describe how to 
construct a first countable compact space $L_\omega$ such that 
$E_0(L_\omega)$ contains only the constant functions.
Let $D \subseteq \CC$
be the closed unit disk; $D$ will be a subspace of $L_\omega$.
We shall first focus on the easier task of constructing
a space $L_2$ such
that $D \subset L_2$ and each $f \in E_0(L_2)$ is constant on $D$.  
After explaining this, we shall
iterate  the procedure to produce $L_\omega$.

Before we build $L_2$, 
we shall show that every nonconstant function $f \in C(D)$
is $1 - 1$ on ``many'' disjoint Cantor sets.  Then, by gluing new disks
on those Cantor sets to form $L_2$, we can make sure that no such $f$ can
extend to a function in $E_0(L_2)$.

For $\theta \in [0,2\pi)$, let $R_\theta$ denote the {\it ray\/}
$\{z \in D: z \ne 0 \,\&\, \arg(z) = \theta\}$.  Let $\cc = 2^{\aleph_0}$.

\medskip

{\bf 2.1. Lemma.}  If $f \in C(D)$ is nonconstant, then there are
$\cc$ distinct $\theta$ such that $f$ is nonconstant on 
$R_\theta$.

{\bf Proof.} The set of all such $\theta$ is open.  \eop

\medskip

We identify $\cc$ with a von Neumann ordinal, so that we may use
$\cc$ also as an index set.

\medskip

{\bf 2.2. Lemma.}  There is a {\it disjoint} family 
$\{ K_\alpha \subset D \backslash \{0\} : \alpha \in \cc \}$
of $\cc$ Cantor sets, with
the following property:
For each nonconstant $f \in C(D)$, there is
a Cantor set $H_f \subset D \backslash \{0\}$ 
such that $f$ is $1 - 1$ on $H_f$ and such that
$A = \{\alpha\in\cc: K_\alpha \subseteq H_f \}$ has size $\cc$.

{\bf Proof.}  First, applying Lemma 2.1 and transfinite induction,
choose, for each nonconstant $f \in C(D)$, a {\it distinct\/}
$\theta_f \in [0,2\pi)$ such that $f$ is nonconstant on
$R_{\theta_f}$.  Then, applying Lemma 1.1, choose a Cantor set
$H_f \subset R_{\theta_f}$ such that $f$ is $1 - 1$ on
$H_f$.  
Partition each $H_f$ 
into $\cc$ disjoint Cantor sets.  
Since the $H_f$ are all disjoint,  
this gives
us the desired family of $\cc \cdot \cc = \cc$  Cantor sets.  \eop

\medskip

Informally, we now replace each $K_\alpha$ by a copy of $K_\alpha \times D$,
identifying $K_\alpha \times \{0\}$ with the old $K_\alpha$.
For different $\alpha$, we want the $K_\alpha \times D$ to point in
``perpendicular directions''.   To make the notion of ``perpendicular''
formal, we simply embed $L_2$ into a Hilbert space.   Since
we want each ``direction'' to be a whole disk, 
we use a {\it complex\/} Hilbert space to simplify 
the notation.  One could use a real Hilbert space
instead by replacing each unit vector in the following
proof by a pair of unit vectors.  In either case, the following
simple criterion can be used to verify first countability.

\medskip
{\bf 2.3. Lemma.} If $\BB$ is a Hilbert space and $X \subset \BB$
is compact in the weak topology, then $X$ is first countable in the
weak topology iff for each $\vec x \in X$, there is a countable
(or finite) $C_{\vec x} \subset \BB$ such that no
$\vec v \in X \backslash \{\vec x\}$ satisfies
$\forall \vec c \in C_{\vec x} (\vec x \cdot \vec c = \vec v \cdot \vec c)$.

{\bf Proof.} By definition of the weak topology, the stated condition
is equivalent to each $\{\vec x\}$ being a $G_\delta$ set in $X$,
which is equivalent to first countability in a compact space. \eop

\medskip
We remark that the condition of Lemma 2.3 need not imply
first countability when $X$ is not weakly compact.

\medskip

{\bf 2.4. Lemma.} There is a first countable uniform
Eberlein compact space $L_2$
such that $D$ is a retract of $L_2$ and
each $f \in E_0(L_2)$ is constant on $D$.

{\bf Proof.}  Let $\BB$ be a complex Hilbert space with an
orthonormal basis consisting of $\cc$ unit vectors 
$\vec e_\alpha$,  for $\alpha \in \cc$, together
with one more, $\vec e$.
We identify $D$ with its homeomorphic copy,
$D' = \{z \vec e : |z| \le 1\} \subset \BB$.
Let $\pi$ be the perpendicular projection from $\BB$ onto
the one-dimensional subspace spanned by $\vec e$.

Let the $K_\alpha \subset D' \backslash \{\vec0\}$
be as in Lemma 2.2 (replacing the $D$ there by $D'$).
Let $L_2$ be the set of all $\vec x \in \BB$ that satisfy
(1) - (3):

{\parindent = 40 pt
\item{(1)} $|\vec x \cdot \vec e| \le 1$, and, for each $\alpha\in\cc$,
$|\vec x \cdot \vec e_\alpha| \le {1 \over 2}$.

\item{(2)} For all distinct $\alpha, \beta$, either
$\vec x \cdot \vec e_\alpha =  0$ or
$\vec x \cdot \vec e_\beta = 0$.

\item{(3)} For all $\alpha$, either
$\vec x \cdot \vec e_\alpha = 0$ or $\pi(\vec x) \in K_\alpha$.
}

\noindent
So, points of $L_2$ are either of the form $z \vec e$, with $|z| \le 1$,
or of the form $z \vec e + w \vec e_\alpha$, where
$|z| \le 1$,  $|w| \le {1 \over 2}$, and $z \vec e \in K_\alpha$.
In particular, $D' = \pi(L_2) \subset L_2$.

We give $L_2$ the topology inherited from the {\it weak\/} topology on $\BB$.
Note that $L_2$ is weakly closed.  Since $L_2$ is
also norm bounded, $L_2$ is compact.
To see that $L_2$ is first countable,  apply Lemma 2.3.
If $\pi(\vec x)$ is in no $K_\alpha$, set $C_{\vec x} = \{\vec e\}$, while
if $\pi(\vec x)$ is in some $K_\alpha$, this $\alpha$ is unique
(by the disjointness of the $K_\alpha$), and we
set $C_{\vec x} = \{\vec e, \vec e_\alpha \}$.

Let $U_\alpha = L_2 \cap \pi^{-1}(K_\alpha) \backslash K_\alpha = 
\{\vec x \in L_2 : \vec x \cdot \vec e_\alpha \ne 0\}$.  Observe:

{\parindent = 40 pt
\item{i.} $U_\alpha$ is an open subset of $L_2$, but
\item{ii.} For each $\vec x \in K_\alpha$,
$L_2 \cap \pi^{-1}(\{\vec x\})$ is nowhere dense in $L_2$.
}

Now, suppose $f \in E_0(L_2)$.  We show that $f$ is constant on $D'$.
If not, fix a Cantor set $H \subseteq D'$
such that $f$ is $1 - 1$ on $H$ and such that
$A = \{\alpha\in\cc: K_\alpha \subseteq H\}$ has size $\cc$.
Since $f \in E_0(L_2)$, we may, for each $\alpha \in A$, choose a nonempty
open $W_\alpha \subseteq U_\alpha$ such that $f$ is constant on
$W_\alpha$.  Then, applying (ii) above, 
choose two distinct points $\vec x_\alpha$ and $\vec y_\alpha$
in $W_\alpha$ such that $\pi(\vec x_\alpha) \ne \pi(\vec y_\alpha)$.

For each $\alpha \in A$, $(\pi(\vec x_\alpha), \pi(\vec y_\alpha))$
is a point in $\{(\vec v, \vec w) \in H \times H :
\vec v \ne \vec w\}$, which is a second countable space. 
Since $A$ is uncountable, these points have a limit
point in the same space, so we may fix {\it distinct\/}
$\vec v, \vec w \in H$ and a sequence of {\it distinct\/}
elements $\alpha_n$ in $A$ ($n \in\omega$) such that the
$\pi(\vec x_{\alpha_n})$ converge to $\vec v$ and
the $\pi(\vec y_{\alpha_n})$ converge to $\vec w$.
Hence, in the {\it weak\/} topology of $\BB$ and $L_2$, the
$\vec x_{\alpha_n}$ converge to $\vec v$ and
the $\vec y_{\alpha_n}$ converge to $\vec w$.
Since $f(\vec x_{\alpha_n}) = f(\vec y_{\alpha_n}) $, we have
$f(\vec v)  = f(\vec w)$, contradicting that $ f$ was $1 - 1$ on $H$.  \eop

\medskip

A similar use of Cantor sets occurs in the construction
in [\RuRu], with the following differences:  Their
$K_\alpha$ were not disjoint; in fact, in [\RuRu]
it appears necessary that every Cantor set gets listed uncountably
many times.
As a result, the space constructed was not first countable.
However, if one does not care about disjointness, there is no
advantage to using a disk, so [\RuRu] used an interval where we
used $D$.  The extra dimension in $D$ lets us prove Lemma 2.2,
which is easily seen to be false of $[0,1]$.
Actually, when the $K_\alpha$ are disjoint, condition (2) above is
redundant, since it follows from (3), but if the $K_\alpha$
are not disjoint, (2) is required to guarantee that $L_2$
is norm bounded.

By iterating our construction, we now prove the following theorem.

\medskip

{\bf 2.5. Theorem.}  There is a first countable uniform Eberlein
compact space
$L_\omega$ such that every function in $E_0(L_\omega)$ is constant.

\medskip

Observe that this is not true for the $L_2$  of Lemma 2.4.
For example, let $g \in E_0(D)$ be nonconstant,
and define $f$ by $f(\vec x) = g(\vec x \cdot \vec e _ \alpha)$.
Then $f \in E_0(L_2)$, and is not constant on $U_\alpha$.
To prevent such functions from
existing, we shall, for each $\alpha$\/: take disjoint Cantor sets
$K_{\alpha\beta} \subset U_\alpha$, and,
for each $\beta$,
attach a new disk going off in a new direction,
labeled by a unit vector $\vec e _{\alpha\beta}$. 
This would create a space $L_3$.  But now, we must iterate this
procedure, to take care of functions on these new disks.
Iterating $\omega$ times, we have
unit vectors $\vec e _t$ indexed by finite sequences from $\cc$.

To describe $L_\omega$, we use the same tree notation 
as in \S1, where now $\cc$ is our index set.
For the rest of this section, let $\BB$ be a complex Hilbert
space with an orthonormal basis consisting of unit vectors
$\{\vec e _ s : s \in \cclw\}$.  We shall use $\vec e$ to abbreviate
$\vec e_{()}$ and $\vec e_\alpha$ to abbreviate $\vec e_{(\alpha)}$.
Let $\pi_n$ be the perpendicular projection from $\BB$ onto
the subspace spanned by $\{\vec e_s: lh(s) < n\}$.  In particular,
$\pi_0(\vec x) = \vec 0$ for all $\vec x$, and $\pi_1$ is the
projection onto the one-dimensional subspace spanned by $\vec e$.

If $lh(s) = n$, let $D_s$ be the set of vectors of the
form $\sum_{i \le n} z_i \vec e _ {s \rest i}$, where each $|z_i| \le 2^{-i}$.
Since $D_s$ is finite dimensional, the weak and norm topologies
agree on $D_s$, and $D_s$ is homeomorphic to $D^{n+1}$.  In particular,
$D_{()} = \{z \vec e: |z| \le 1\}$ plays the role of the
$D'$ in the proof of Lemma 2.4.
Note that if $i \le n$, then $\pi_{i+1}(D_s) = D_{s \rest i}$.

We begin by enumerating enough of the conditions required of
the Cantor sets $K_t$ ($t \in \cclw$) to define $L_\omega$. 
Then, after defining
$L_\omega$, we prove a sequence of lemmas, adding conditions on 
the $K_t$ as necessary, to show $L_\omega$ has the desired properties.

\medskip

{\bf 2.6. Basic requirements on the $K_t$.} 
\itemitem{\bf (Ra)}  $K_{()} = \{\vec0\}$.
\itemitem{\bf (Rb)}  For each $s$,  the $K_{s\alpha}$ for $\alpha\in\cc$ are
disjoint closed subsets of $D_s$, and $\vec x \cdot \vec e_s \ne 0$
for all $\vec x \in K_{s\alpha}$.
\itemitem{\bf (Rc)}  For each $s$ and each $\beta$,  if $n = lh(s)$,
then $\pi_n(K_{s\beta}) \subseteq K_s$.

\medskip

In particular, for $s = ()$, we have $K_{\alpha} \subset D_{()}$,
as in the proof of Lemma 2.4.  Now, we iterate that construction
by using the $K_{\alpha\beta}$, $K_{\alpha\beta\gamma}$, etc.
The $K_{()} = \{\vec0\}$ plays no role in the definition of $L_\omega$, but
is included to make some of the notation more uniform.
Item (Rc) for $n = 0$ says nothing; for $n = 1$,
$\pi_1(K_{\alpha\beta}) \subseteq K_{\alpha}$ corresponds
to the informal idea above that the $K_{\alpha\beta}$ are chosen
inside $U_\alpha$.

We shall need to add conditions (Rd)(Re) to (Ra)(Rb)(Rc) later.

\medskip
\vbox{  
{\bf 2.7 Definition.}  $L_\omega$ is the set of all $\vec x \in \BB$ 
that satisfy (1) - (3):

\item{(1)} For each $s$, $|\vec x \cdot \vec e_s| \le 2^{-lh(s)}$.

\item{(2)} For all $s,t$ such that $s\perp t$,
$\vec x \cdot \vec e_s =  0$ or
$\vec x \cdot \vec e_t = 0$.

\item{(3)} For all $t$,  if $n = lh(t)$, then either
$\vec x \cdot \vec e_{t} = 0$ or
$\pi_n(\vec x) \in K_t$.

\noindent
We give $L_\omega$ the weak topology.  $L_n = \pi_n(L_\omega)$.
For $\vec x \in L_\omega$,
$P(\vec x) = \{s \in \cclw : \vec x \cdot \vec e_s \ne  0\}$.
For $t \in \cclw$ and $n = lh(t)$, set
$U_t = L_\omega \cap (\pi_n^{-1}(K_t)\backslash K_t)$.
}

\medskip

{\bf 2.8. Lemma.}  Each $L_n$ is a closed subset of $L_\omega$ and
$\bigcup_{n\in\omega}L_n$ is dense in $L_\omega$.

{\bf Proof.} $L_n \subseteq L_\omega$ holds 
because each of (1), (2), (3) is preserved under $\pi_n$.
Density follows because for every
$\vec x \in \BB$, the $\pi_n(\vec x)$ converge weakly
(and in norm) to $\vec x$.
$L_n$ is closed in $L_\omega$ because $\pi_n(\BB)$ is weakly 
closed in $\BB$.  \eop

\medskip

We think of the $L_n$ as the {\it levels\/} in the construction.
$L_0 = K_{()}$.  $L_1 = D_{()}$.  $L_2$ is exactly the
space constructed in the proof of Lemma 2.4.
The $U_t$ will play the same role here as the $U_\alpha$ did there.
Elements of $L_3 \backslash L_2$ are of the form
$r_0 \vec e + r_1 \vec e_\alpha + r_2 \vec e_{\alpha\beta}$, where
$0 < |r_i| \le 2^{-i}$ for each $i$,
$r_0 \vec e \in K_\alpha$, and
$r_0 \vec e + r_1 \vec e_\alpha \in K_{\alpha\beta}$.

\medskip

{\bf 2.9. Lemma.} 
\item{i.}  For each $\vec x \in L_\omega$, $P(\vec x)$ is a path in $\cclw$.
\item{ii.}  For each $\vec x \in L_\omega$, $\|\vec x\|^2 \le {4 \over 3}$.
\item{iii.}  $L_\omega$ is weakly closed in $\BB$.
\item{iv.}  $L_\omega$ is first countable and compact.
\item{v.}  Each $U_t$ is open in $L_\omega$.

{\bf Proof} For (i), use items (2),(3) in the definition of $L_\omega$ and
the fact that $\vec x \cdot \vec e_s \ne 0$ for all $\vec x \in K_{s\alpha}$.
Now, (ii) follows by item (1).
(iii) is immediate from the definition of $L_\omega$, and compactness
of $L_\omega$ follows by (iii) and (ii).  First countability
follows from Lemma 2.3; $C_{\vec x} = \{ \vec e_s : s \in P(\vec x) \}$, unless 
$P(\vec x)$ is finite with maximal element $s$ and $\vec x \in K_{s\alpha}$,
in which case 
$C_{\vec x} = \{ \vec e_s : s \in P(\vec x) \} \cup \{\vec e _ {s\alpha}\}$.
For (v), note
that $U_t = \{\vec x \in L_\omega : \vec x \cdot \vec e_t \ne 0\}$.  \eop

\medskip

Applying conditions (Rc) and (Rb) on the $K_s$, we have the
following lemma.

\medskip

{\bf 2.10. Lemma.} 
\item{i.}  For each $t$, if $n \le lh(t)$ and $s = t \rest n$,
then $K_s \supseteq \pi_n(K_t)$.
\item{ii.}  Each $K_t \subseteq L_{lh(t)}$.

\medskip

If the $K_\alpha$ are chosen as in the proof of Lemma 2.4,
then every $f \in E_0(L_2)$ will be constant on $D_{()}$.  We must
be careful not to destroy this property in choosing the $K_{\alpha\beta}$
and passing to $L_3$.  In the proof of Lemma 2.4, it was important that each
$\pi^{-1}(\{\vec x\})$ was nowhere dense.  Now, $L_2 \cap \pi^{-1}(\{\vec x\})$
will still be nowhere dense in $L_2$, but depending on how the
$K_{\alpha\beta}$ meet this set, $L_3 \cap \pi^{-1}(\{\vec x\})$
might have interior points.  To handle this, we assume the
following product structure on the $K_s$:

\medskip

\itemitem{\bf (Rd)}  For each $s$ of
length $n \ge 0$ and each  $\alpha$, there are
a nonempty relatively
clopen subset $P \subseteq K_s$ and a homeomorphism $\psi$
from $P \times 2^{\omega}$ onto $K_{s\alpha}$,  satisfying
$\pi_n(\psi(\vec x, y)) = \vec x$  for all $\vec x \in P$
and all $y \in 2^\omega$.

\medskip

\noindent
Note that (Rd) implies that $\pi_n(K_{s\alpha}) = P$.
Induction on $lh(s)$ establishes the next lemma.

\medskip
{\bf 2.11. Lemma.}  $K_s$ is a Cantor set whenever $lh(s) > 0$.

\medskip
{\bf 2.12. Lemma.}  Suppose that $m > 0$ and $C$ is a closed subset
of $L_m$ such that $C$ is nowhere dense (in the relative topology
of $L_m$) and $C \cap K_s$ is nowhere dense (in the relative
topology) in $K_s$ for all $s$ of length $m$.  Then
$L_\omega \cap \pi_m^{-1}(C)$ is nowhere dense in $L_\omega$.  In particular, 
$L_\omega \cap \pi_m^{-1}(\{\vec x\})$ is nowhere
dense in $L_\omega$ for all $\vec x \in L_m$.

{\bf Proof.}  The ``in particular'' follows from Lemma 2.11, 
which implies that $C = \{\vec x\}$ satisfies the hypotheses of Lemma 2.12.
Now set $C_n = L_n \cap \pi_m^{-1}(C)$ for each
$n \ge m$; so $C_m = C$. 
To prove 2.12, 
since $\bigcup_{n\in\omega}L_n$ is dense in $L_\omega$, it
suffices to prove claim (i) below.   
To do this, we prove claims (i) and (ii) together,
by induction on $n \ge m$. 

\medskip
\item{i.} For each $n \ge m$, $C_n$ is nowhere dense in $L_n$.
\item{ii.} Whenever $lh(s) = n$, 
$C_n \cap K_s$ is nowhere dense in $K_s$.

\medskip

\noindent
Claim (ii) for $n+1$ follows from (ii) for $n$ plus assumption (Rd)
on the $K_s$, and claim (i) for $n+1$ follows from (i) and (ii) for $n$
(just using (Ra),(Rb),(Rc)).  \eop

\medskip

For each $s \in \cclw$, with $lh(s) = n$,
let
$$
\hat K_s = \{\vec v + z \vec e_s :
\vec v \in K_s \,\&\, |z| \le 2^{-n} \} \ \ .
$$
Note that $\hat K_s$ is
homeomorphic to $K_s \times D$ and is a subset of $L_\omega$.
If $H \subseteq \hat K_s$ and $\vec v \in K_s$, let
$H_{\vec v}$ be the ``vertical slice'', 
$ \{\vec v + z \vec e_s :  |z| \le 2^{-n} \}$.
Call a function $f$ {\it $s$-level-constant\/}
iff $f$ only depends on the $\vec v$ here;  that is, 
$f$ is constant on each $(\hat K_s)_{\vec v}$.
In particular, $f$ is $()$-level-constant iff $f$ is constant
on $D_{()}$, and the $K_\alpha$ chosen as in the proof of 
Lemma 2.4 will ensure that every $f \in E_0(L_\omega)$ is 
$()$-level-constant.  Likewise, we shall choose the 
$K_{s\alpha}$ to ensure that every $f \in E_0(L_\omega)$ is 
$s$-level-constant.  
Note first that if we do this for all $s$, then $f$ is constant.

\medskip
{\bf 2.13. Lemma.}  If $f \in C(L_\omega)$ is $s$-level-constant for all
$s \in \cclw$, then $f$ is constant.

{\bf Proof.} By induction on $n$, $f$ is constant on each $L_n$.
The result follows because 
$\bigcup_{n\in\omega}L_n$ is dense in $L_\omega$.  \eop

\medskip

Now we list the final condition on the $K_{s\alpha}$:

\medskip

\itemitem{\bf (Re)}  For each $s$ of length $n$ and each $f \in C(L_\omega)$:
If $f$ is not $s$-level-constant, then there are a nonempty clopen set
$P \subseteq K_s$, a Cantor set 
$H \subseteq \{\vec v + z \vec e_s: \vec v \in P \,\&\, |z| \le 2^{-n}\}$,
and uncountably many different $\alpha$
such that $K_{s\alpha} \subset H$,  and
for each $\vec v \in P$, 
$f$ is $1 - 1$ on $H_{\vec v}$.  

\medskip

We must verify that we may choose 
the $K_t$ to meet all five conditions (Ra), (Rb), (Rc), (Rd), (Re).
We choose these by induction on $lh(t)$.  Condition (Ra) specifies
$K_{()}$, and the $K_\alpha$ will be exactly as in the proof
of Lemma 2.4; these were chosen by applying Lemma 2.2.  Likewise,
given $K_s$ with $lh(s) > 0$, we choose the $K_{s\alpha}$
by applying the next lemma to $\hat K_s$.
In fact, we modify the proof of Lemma 2.2, 
replacing Lemma 1.1 by Lemma 1.2, to prove this lemma.

\medskip

{\bf 2.14. Lemma.}  
Let $\{E_\delta : \delta\in\cc \}$ be
a partition of $2^\omega$ into $\cc$ Cantor sets.
If $K$ is a Cantor set, then  there 
is a disjoint family 
$\{ K_\alpha \subset K \times (D\backslash\{0\}) : \alpha\in\cc \}$ 
of $\cc$ Cantor sets 
with the following property:  For each $f \in C(K \times D)$
with $f \rest (\{x\} \times D)$ nonconstant for {\it some} 
$x \in K$, there are a nonempty clopen $P \subseteq K$ and
an $H \subset P \times (D\backslash\{0\})$ that satisfy
conditions (1) - (4):

\item{(1)} $H_x$ is a Cantor set for every $x \in P$.
\item{(2)} $f$ is $1 - 1$ on $\{x\} \times H_x$ for every $x \in P$.
\item{(3)} There is a continuous $\varphi : H \to 2^\omega$
such that the map $(x,y) \mapsto (x, \varphi(x,y))$ is a homeomorphism
from $H$ onto $P \times 2^\omega$.
\item{(4)} For each $\delta\in\cc$, the set
$\{(x,y) \in H: \varphi(x,y) \in E_\delta\}$ is one of the $K_\alpha$.

{\bf Proof.}  
First, for each such $f$, apply continuity
to choose a nonempty clopen $P_f \subseteq K$ such that 
for $\cc$ different $\theta \in [0,2\pi)$, 
$f \rest (\{x\} \times R_\theta)$ fails to be constant for all
$x \in P_f$.  Then, by transfinite induction, choose a distinct
$\theta_f$ for each such $f$ such that 
$f \rest (\{x\} \times R_{\theta_f})$ is not constant for all $x \in P_f$.  
Then, choose $H_f \subseteq P_f \times R_{\theta_f}$ such that
(1), (2), and (3) hold; this is possible by Lemma 1.2.
Of course, $\varphi = \varphi_f$ depends on $f$.
Finally, let the $K_\alpha$ enumerate all
the sets $\{(x,y) \in H_f: \varphi_f(x,y) \in E_\delta\}$ as $f$
and $\delta$ vary. \eop

\medskip

Now we complete the proof of Theorem 2.5.

\medskip

{\bf Proof of Theorem 2.5.}  
Construct $L_\omega$ as above.  Suppose $f \in E_0(L_\omega)$.
By Lemma 2.13, it suffices to prove that $f$ is
$s$-level-constant for each $s$.  Suppose not.  Fix $H,P$
as in condition (Re) above, so that $A = \{\alpha : K_{s\alpha} \subseteq H\}$
is uncountable.  For $\alpha \in A$, choose a nonempty open $W_\alpha$
such that $W_\alpha \subseteq \overline W_\alpha \subseteq U_{s\alpha}$ and
$f$ is constant on $\overline W_\alpha$.  Then
$\pi_{n+1}(\overline W_\alpha) \subseteq K_{s\alpha} \subseteq H$
and $\pi_n(\overline W_\alpha) \subseteq \pi_n(K_{s\alpha}) \subseteq P
\subseteq K_s \subseteq L_n$.  Choose $\vec x_\alpha$
and $\vec y_\alpha$ in $\overline W_\alpha$ such that
$\pi_n(\vec x_\alpha) = \pi_n(\vec y_\alpha)$ but
$\pi_{n+1}(\vec x_\alpha) \ne \pi_{n+1}(\vec y_\alpha)$;
this is possible because $\pi_{n+1}(\overline W_\alpha)$ is closed in 
$K_{s\alpha}$ and, by Lemma 2.12, is not nowhere dense in
$K_{s\alpha}$.  As in the proof of Lemma 2.4,
there are distinct $\vec v, \vec w \in H$ and a
sequence of distinct
elements $\alpha_k$ in $A$ ($k \in\omega$) such that the
$\pi_{n+1}(\vec x_{\alpha_k})$ converge to $\vec v$ and
the $\pi_{n+1}(\vec y_{\alpha_k})$ converge to $\vec w$.
Then, in the weak topology, the
$\vec x_{\alpha_k}$ converge to $\vec v$ and
the $\vec y_{\alpha_k}$ converge to $\vec w$.
So, $f(\vec v) = f(\vec w)$, while 
$\pi_n(\vec v) = \pi_n(\vec w) \in \pi_n(H) \subseteq P$,
contradicting that $f$ is $1 - 1$ on $H_{\pi_n(\vec v)}$.  \eop

\medskip

Finally, we remark on $E_0(X,M)$ for other $M$.

\medskip

{\bf 2.15. Lemma.}
If $X$ is a compact Hausdorff space and $M$ is any Hausdorff space,
then
\item\item{(1)} $E_0(X,\RR)$ contains only the constant functions 

\noindent
implies
\item\item{(2)} $E_0(X,M)$ contains only the constant functions.

\noindent
If $M$ contains a closed interval, then (2) implies (1).

{\bf Proof.} For $(1) \to (2)$, fix $f \in E_0(X,M)$.  We may
assume $M = f(X)$, whence $M$ is compact.  For each $g: M \to [0,1]$,
$g \circ f$ is in $E_0(X,\RR)$ and hence constant, which implies that $f$
is constant.  For $\neg (1) \to \neg (2)$, if $g$ maps $\RR$ homeomorphically
into $M$ and $f$ is a nonconstant function in $E_0(X,\RR)$, then
$g \circ f$ is a nonconstant function in $E_0(X,M)$. \eop

\medskip

In particular, in making $E(X)=
E_0(X,\RR)$ small, we also make $E_0(X,\CC)$ small.   
Note that 2.15 can fail if $M$ does
not contain an interval, since then, if $X$ is a closed interval,
$E_0(X,M) = C(X,M)$ contains only the constant functions
(since every arc contains a simple arc),
while $E_0(X,\RR)$ is dense in $C(X,\RR)$.  We do not 
study 2.15 for such $M$ in detail here, but it seems to involve the
geometric-topological properties of $X$ and $M$.

\bigskip

{\bf \S3. Making $E_0(X)$ Big.}  Here, we consider spaces $X$ for
which $E_0(X,M) = C(X,M)$.  This turns out to be an interesting
topological property of $X$.  We begin with a simple remark.

The condition $E_0(X,M) = C(X,M)$ is not hereditary to closed
subsets of $X$, but it is, in many cases, hereditary to 
{\it regular closed\/} subsets -- that is, to subsets of the
form $\overline U$, where $U$ is open in $X$.

\medskip

{\bf 3.1. Lemma.} Suppose that $X$ is a compact Hausdorff space,
$E_0(X,\RR) = C(X,\RR)$, and $Y$ is a regular
closed subspace of $X$.  Then
$E_0(Y,\RR) = C(Y,\RR)$.

{\bf Proof.}
Say $Y = \cl U$, where $U$ is open.  Suppose $g \in C(Y,\RR)$.
By the Tietze Extension Theorem, $g$ can be extended to an
$f  \in C(X,\RR)$.  Then $\Omega_g \cap U = \Omega_f \cap U$.
Since $E_0(X,\RR) = C(X,\RR)$, we have that $\Omega_f$ is dense in $X$,
so $\Omega_g$ is dense in $Y$.  \eop

\medskip

We remark that in Lemma 3.1, one can replace $\RR$ by any Banach
space (using a slightly longer proof), but not by an arbitrary metric
space $M$.  For a counter-example, let $M$ be a Cantor set
and  let $X$ be the cone over $M$.  Then $E_0(X,M) = C(X,M)$ contains
only constant functions.  But $X$ contains a regular closed
$Y$ homeomorphic to $M \times [0,1]$, and
$E_0(Y,M) \ne C(Y,M)$.  Also, even in the simple case
$M = \RR$, the property $E_0(Y,\RR) = C(Y,\RR)$ holds
for {\it all\/} closed $Y \subseteq X$ iff $X$ is scattered;
if $X$ is not scattered, then $X$ will
contain a closed subset $Y$ which
is separable with no isolated points, which implies $E_0(Y,\RR) \ne C(Y,\RR)$
(by $(2) \imp (1)$ of Theorem 3.2 below).

Now, to study the property $E_0(X,M) = C(X,M)$, 
it is convenient to generalize our notions in two ways.

First, although $X$ will always be compact and $M$ will always
be metric, we look at more general functions from $X$ into $M$.
In particular recall that $f : X \to M$ is called 
{\it Borel measurable\/} iff the inverse image of every open set
is a Borel subset of $X$, and 
{\it Baire measurable\/} iff the inverse image of every open set
is a Baire subset of $X$; the Baire sets are the 
$\sigma$-algebra generated by the open $F_\sigma$ sets.
The Baire measurable functions into a separable Banach space
form the least class of functions containing the continuous
functions and closed under pointwise limits.

Second, we consider also
$\widehat\Omega_f$, which we define to be the union of
all open $U \subseteq X$ such that for some first category 
set $C \subseteq X$, 
$f$ is constant on $U\backslash C$.
Note that regardless of $f$, $\Omega_f$ (defined in the Introduction)
and $\widehat\Omega_f$ are open, with $\Omega_f \subseteq \widehat\Omega_f$.
If $f$ is continuous, then $\Omega_f = \widehat\Omega_f$.

The property $E_0(X,\RR) = C(X,\RR)$ is just one of a sequence
of related properties:

\medskip

\vbox{
\item{\bf (1)}  Every nonempty open subset of $X$ is either nonseparable
or contains an isolated point.
\item{\bf (2)}  $E_0(X,\RR) = C(X,\RR)$.
\item{$\bf (3)$}  For all metric spaces $M$, $E_0(X,M) = C(X,M)$.
\item{$\bf (3')$}  For all separable metric spaces $M$ and all
Baire measurable $f : X \to M$, $\widehat \Omega_f$ is dense in $X$.
\item{$\bf (3'')$}  For all separable metric spaces $M$ and all
Baire measurable $f : X \to M$, $\Omega_f$ is dense in $X$.
\item{\bf (4)}  For all separable metric spaces $M$ and all
Borel measurable $f : X \to M$, $\widehat \Omega_f$ is dense in $X$.
\item{\bf (5)}  For all separable metric spaces $M$ and all
Borel measurable $f : X \to M$, $\Omega_f$ is dense in $X$.
\item{\bf (6)}  In $X$, every nonempty $G_\delta$ set
has a  nonempty interior.
\smallskip
\item{$\bf (*)$} In $X$, every first category set is nowhere dense.
}

\medskip

Conditions (1) -- (6) are listed in order of increasing strength.
Condition $(*)$ does not fit into the sequence, but is relevant
by the next Theorem.

\medskip

{\bf 3.2. Theorem.} Suppose $X$ is compact Hausdorff.  Then
$$
(6) \imp (5) \imp (4)  \imp (3)  \iff (3')  \iff (3'')  \imp (2)  \imp (1) \ \ .
$$
Furthermore (5) is equivalent to $(*)$ plus (4).

{\bf Proof.} For $(2) \imp (1)$, assume (1) fails; so there is nonempty open
$U$ which is separable and has no isolated points.  Let 
$Y = \overline U$.  By Lemma 3.1, it is sufficient to produce an
$f \in  C(Y,\RR) \setminus E_0(Y,\RR)$.
Let $\{p_n : n \in \omega\}$ be dense in $U$, hence in $Y$.
For each distinct $m,n$, $\{f \in C(Y,\RR) : f(p_m) \ne f(p_n) \}$ is 
dense and open in $C(Y,\RR)$ (in the usual norm topology), so by
the Baire Category Theorem, there is an $f \in  C(Y,\RR)$ such that
$f(p_m) \ne f(p_n)$ whenever $m \ne n$.  But then for each 
$r \in \RR$, $f^{-1}\{r\}$ contains at most one $p_n$, and is hence
nowhere dense in $Y$, since $Y$ has no isolated points.
Thus, $f \notin E_0(Y,\RR)$.    

Clearly, $(3'')  \imp (3')  \imp (3)$, so to prove these three
are equivalent, we assume (3), fix a Baire measurable $f: X \to M$, and 
show that $\Omega_f$ is dense in $X$.  
Since $M$ can be embedded into a separable Banach space, we may
assume that $M$ {\it is\/} a Banach space; now, we can
let $g_n : X \to M$, for $n \in \omega$,
be continuous functions such that $f$ can be obtained from
the $g_n$ by some transfinite iteration of taking pointwise limits.
Define $g : X \to M^\omega$ by: $g(x)_n = g_n(x)$.  Then
$\Omega_g \subseteq \Omega_f$, and, by (3), 
$\Omega_g$ is dense in $X$.

To prove $(6) \imp (5)$, observe that for {\it any\/} compact $X$,
if $H$ is a nonempty closed $G_\delta$ and $f$ is a Borel measurable
map into a second countable space, there
is always a nonempty closed $G_\delta$ set $K \subseteq H$
such that $f$ is constant on $K$.

The rest of the chain of implications from (6) down
to (1) are now trivial.  To see that $(5) \imp (*)$, let $C$
be first category; then $C \subseteq \bigcup_{n \in \omega} K_n$,
where each $K_n$ is closed nowhere dense.  Define $f: X \to 2^\omega$
so that $f(x)_n$ is 1 if $x \in K_n$ and 0 if $x \notin K_n$.
Then $\Omega_f$ is dense and open, and is disjoint from all the $K_n$,
so $C$ is nowhere dense.

To see that $(*)$ plus (4) implies (5), we let  $f$ be Borel measurable;
to prove $\Omega_f$ dense, we fix a nonempty open $V$ and try to
find a nonempty open $U \subseteq V$ such that $f$ is constant on $U$.
By (4), there is a nonempty open $W \subseteq V$ such that 
$f$ is constant on $W \backslash C$ for some first category $C$.
By $(*)$, $C$ is nowhere dense, so let $U = W \backslash \overline C$.  \eop

\medskip

A familiar example of a space satisfying (6) is $\beta \NN \backslash \NN$.

Conditions (5), (4), $(3'')$, and $(3')$ involve arbitrary
Baire or Borel measurable maps.  Each of these conditions
is equivalent to the restatement we obtain by replacing 
$M$ by the Cantor set $2^\omega$.  This is easily seen by translating
the condition to one involving an $\omega$-sequence of Borel or Baire
sets.  For example, (5) is equivalent to the statement that given
Borel sets $B_n$ ($n \in  \omega$), the union of all open
$U$ such that $\forall n (U \subseteq B_n \ or \ U \cap B_n = \emptyset)$
is dense in $X$.

This is not true for (3), which involves continuous functions.
For example, if $X$ is connected, then,
trivially, $E_0(X,2^\omega) = C(X,2^\omega)$, whereas 
$E_0(X,\RR)$ need not be all of $C(X,\RR)$.
If $X$ is zero-dimensional, then $E_0(X,2^\omega) = C(X,2^\omega)$ does
imply (3).  In fact, for zero-dimensional spaces, (3) has a
restatement in terms of sequences of clopen sets (see the proof
of Theorem 3.3(c) below).

Regarding $(3) \imp (2)$, if $E_0(X,M) = C(X,M)$ for
{\it any\/} $M$ containing an interval, then (2) holds.
In \S4, we show that (2) does not imply (3), 
although it is easy to see that (2) implies
that $E_0(X,\RR^n) = C(X,\RR^n)$ for each finite $n$.
Counter-examples to the other implications of Theorem 3.2
reversing are provided by some
fairly familiar spaces, as we point out below.
However, the implications do reverse for certain families of spaces.
In particular, we consider the cases when $X$ is
{\it extremally disconnected\/} ({\it e.d.}),
when $X$ is an {\it Eberlein compact},
when $X$ is a {\it LOTS},
and when has the {\it ccc}.
$X$ is called e.d. iff
the closure of every open subset of $X$ is clopen.
$X$ is an Eberlein compact iff $X$ is homeomorphic to a weakly
compact subspace of a Banach space.
$X$ is a LOTS iff $X$ is a totally ordered set, given the order topology.
$X$ has the ccc iff there is no uncountable family of disjoint
open sets in $X$.

The following theorem summarizes what we know for these and some
other simple classes.

\medskip

{\bf 3.3.  Theorem.}  Let $X$ be compact Hausdorff.
\item{a.} If $X$ is metric, then $(1) \iff (5)$, and (1) -- (5) hold
iff the isolated points of $X$ are dense in $X$.
\item{b.} If $X$ is e.d., then $(2) \iff (4)$.
\item{c.} If $X$ is zero-dimensional, then $(2) \iff (3)$.
\item{d.} If $X$ is ccc, then $(3) \iff (5)$.
\item{e.} If $X$ is Eberlein compact, then $(4) \iff (5)$, and (4) -- (5) hold
iff the isolated points of $X$ are dense in $X$.
\item{f.} If $X$ is a LOTS, then $(1) \iff (3)$.

{\bf Proof.} (a) is immediate from the fact that compact metric
spaces are separable.

For (b), assume (2), and let $f : X \to M$ be Borel measurable.
Let $\{B_n : n \in \omega\}$ be an open base for $M$.
Since each $f^{-1}(B_n)$ is a Borel set, there are open
$U_i \subseteq X$, for $i \in \omega$, such that each $f^{-1}(B_n)$
is in the $\sigma$-algebra generated by 
$\{U_i : i \in \omega\}$.  Let $K_i = \overline U_i$, which is clopen.
Define $g : X \to 2^\omega$ so that
$g(x)_i = 1$ iff $x \in K_i$.  Since $2^\omega$ is embeddable in $\RR$,
(2) implies that $\Omega_g$ is dense.   Since
$\bigcup_{i \in \omega}(K_i\backslash U_i)$ is first category,
$\Omega_g \subseteq \widehat\Omega_f$, so
$\widehat\Omega_f$ is dense.

For (c), assume (2), and let $f : X \to M$ be continuous.
Let the $B_n$ be as in the proof of (b).
Since each $f^{-1}(B_n)$ is an open $F_\sigma$ set, there are
clopen sets $K_i \subseteq X$ for $i \in \omega$ such that each
$f^{-1}(B_n)$ is a union of some subfamily of the $K_i$.  Now,
construct $g$ as in the proof of (b), and note that
$\Omega_g \subseteq \Omega_f$.

For (d), assume $(3'')$, and let $f : X \to M$ be 
Borel measurable.  Since $X$ is ccc, there is a Baire measurable
$g : X \to M$  and a Baire first category set $C$ such that
$f(x) = g(x)$ for all $x \notin C$.  Define
$h : X \to M \times \{0,1\}$ so that 
$h(x) = (g(x), 0)$ if $x \notin C$, and $h(x) = (g(x), 1)$ if $x \in C$.
Then, applying $(3'')$ to $h$, $\Omega_h$ is dense in $X$.  Since
$\Omega_h \subseteq \Omega_g$ and
$\Omega_h \cap C = \emptyset$,
$\Omega_f$ is dense in $X$.

For (e), assume that $X$ is Eberlein compact
and satisfies (4); we prove that the isolated points are dense. 
By a result of Benyamini, Rudin, and Wage [\BRW],
there is a dense $G_\delta$ set $Y \subseteq X$ such that $Y$ is
metrizable in its relative topology.  Fix some metric on $Y$; then for
$E \subseteq Y$, $diam(E)$ denotes the diameter of $E$ 
with respect to this metric.
For each $n$, let $\WW_n$ be a maximal disjoint family of open nonempty
subsets of $Y$ of diameter $\le 2^{-n}$; 
then $W_n = \bigcup\{W: W \in \WW_n\}$ is
open and dense.  Assume also that each $\WW_{n+1}$ refines $\WW_n$ in
the sense that $\forall W \in \WW_{n+1} \exists V \in \WW_n (W \subseteq V)$,
and for each $V \in \WW_n$ which is not a singleton, there are at least
two $W \in \WW_{n+1}$ such that $W \subseteq V$.  Let
$Z = \bigcap_n W_n$; then $Z$ is also a dense $G_\delta$ subset of $X$.
For each $n$, let $f_n: Z \to 2$ be any function such that 
$f_n$ is constant on {\it every\/} $W \in \WW_{n+1}$ and 
$f_n$ is constant on {\it no\/} $V \in \WW_{n}$ unless $V$ is a singleton.
This defines $f: Z \to 2^\omega$ by $f(z)_n = f_n(z)$.
Let $M$ be the disjoint sum of $2^\omega$ and a single point, $p$,
and extend $f$ to a function $\tilde f: X \to M$ 
by mapping $X \backslash Z$ to $p$.
Then $\tilde f$ is Borel measurable, and every point in
$\widehat \Omega_{\tilde f}$ is isolated in $X$.

For (f),  assume (1), and fix $f \in C(X,M)$; we must show
that $\Omega_f$ is dense.  So, fix a
nonempty open interval $(a,b) \subseteq X$.  We must produce
a nonempty open $W \subseteq (a,b)$ such that $f$ is constant on 
$W$.  This is trivial if $(a,b)$ contains an isolated point, so
assume that $(a,b)$ contains no isolated points, and hence is
nonseparable.  For each $n$, there is a finite cover of
$[a,b]$ by open intervals, $I^n_1, I^n_2, \ldots$ such that
each $diam(f(I^n_j)) \le 1/n$.  Since $(a,b)$ is
nonseparable, we can choose $W \subseteq (a,b)$ to be an open
interval which contains none of the endpoints of any $I^n_j$.
Then for each $n$, $W$ is a subset of some $I^n_j$, so
$diam(f(W)) \le 1/n$.  Thus, $f$ is constant on $W$. \eop

\medskip

A (compact) Suslin line in which every open interval is nonseparable
is a ccc LOTS which satisfies (1), and hence (5),
applying (d) and (f) of the Theorem.  Of course, the Suslin line
does not satisfy (6).  The absolute (or projective cover) of
a Suslin line is a compact ccc e.d. space which satisfies (5) but
not (6).   So is $\beta\NN$, but this example
is ``trivial'' because the isolated points are dense.  Note, however,
that it is consistent with the axioms of set theory that
there are no Suslin lines, in which case (5) 
for a ccc space would imply that the isolated points are dense.

In general, for a LOTS,
(1) need not imply (4).  A simple counter-example is
$X = [0,1]^\omega$, ordered lexically; (4) is refuted by
$f(x) = \sum_{n\in\omega}x_n \cdot 2^{-n}$.
One can replace $[0,1]$ by the Cantor set here to get a
zero-dimensional LOTS, providing also a counter-example to
(c) extending to $(2) \iff (4)$.

The Stone space of an atomless probability algebra 
is a compact e.d. space which satisfies (1) but not (3).
To refute (3), let the $K_i$ ($i \in \omega$) be
clopen independent events of probability ${1 \over 2}$, and
construct $g$ as in the proof of (b) above.
This provides a counter-example to replacing (2) by (1) in 
either (b) or (c).

Conditions (4), $(*)$, and (5) are equivalent to algebraic conditions on
the Boolean algebra of regular open subsets of $X$ (see [\Hal, \Sik]);
in particular, each condition holds for $X$ iff it holds
for the absolute of $X$.
Condition (4) is equivalent to the $(\omega,\omega)$ - distributive
law:  
$$
\bigwedge_{n\in\omega}\bigvee_{i\in\omega}b_{n,i} =
\bigvee\bigl\{\bigwedge_{n\in\omega}b_{n,\varphi(n)} : 
\varphi \in \omega^\omega\bigr\} \ \ .
$$
Condition $(*)$ is equivalent to the weak $(\omega,\infty)$ - distributive
law; that is, for each cardinal $\kappa$,
$$
\bigwedge_{n\in\omega}\bigvee_{\alpha\in\kappa}b_{n,\alpha} =
\bigvee\bigl\{\bigwedge_{n\in\omega}\bigvee_{\alpha\in\varphi(n)} b_{n,\alpha} : 
\varphi \in ([\kappa]^{<\omega})^\omega\bigr\} \ \ .
$$
Here, $[\kappa]^{<\omega}$ is the set of finite subsets of $\kappa$.
Condition (5) is simply (4) plus $(*)$ by Theorem 3.2, which is
equivalent to the standard $(\omega,\infty)$ - distributive law.

Proceeding completely off the deep end, we may regard
the open (or regular open) subsets of $X$
as a forcing order (see a set theory text, such as
[\Jech] or [\Kuna]).  Then (4) is simply the statement that
the order adds no reals, while (5) is the stronger statement that the
order adds no $\omega$ - sequences.
Condition $(*)$  is the finite
approximation property familiar from random real forcing
or Sacks forcing; that is, for each $\kappa$ and each
$\psi: \omega \to \kappa$ in the generic extension,
there is a $\varphi: \omega \to [\kappa]^{<\omega}$ in the ground
model such that each $\psi(n) \in \varphi(n)$.
Prikry forcing at a measurable cardinal (see \S37 of [\Jech])
is an example
of a forcing order (and hence, by the standard translation,
a compact e.d. space) which satisfies (4) but not $(*)$,
and hence not (5).
Another such example is Namba forcing (see \S26 of [\Jech]).

Returning temporarily to Earth, it is natural to ask
which of the properties,  (1) -- (6),  $(*)$, are preserved
by finite products.  Now, (1) and (6) are, trivially.
We don't know about (2), but (3) is; to see this,
identify $C(X\times Y, M)$ with
$C(X, C(Y,M))$, and note that $C(Y,M)$ is another metric space.
Finally, $(5) = (4) + (*)$ is refuted by a well known forcing
order.  Let $S \subset \omega_1$ be stationary and co-stationary.
Let $\PP, \QQ$ be Jensen's forcings for shooting a club through
$S, \omega_1\backslash S$, respectively
(see VII.H25 of [\Kuna]).  Then $\PP,\QQ$ each satisfy (5), while
$\PP \times \QQ$ collapses $\omega_1$, and hence satisfies
neither (4) nor $(*)$.  One may now translate $\PP,\QQ$ into
compact e.d. spaces (by the standard translation), or into
Corson compacta (using the fact that these partial orders
have no decreasing $\omega_1$ chains).

Preservation by infinite products is uninteresting.
If $X$ is an infinite product of spaces with more than one point,
then all of (2) -- (6) fail, as does $(*)$, whereas (1) will hold if,
for example, infinitely many of the $X_n$ are nonseparable.
See Theorem 5.4 for more about such products.

Some of the results in in this section overlap results of Bella, Hager,
Martinez, Woodward, and Zhou [\BHMWZ, \BMW, \MW].
They also defined $E_0(X, \RR)$
(which they called $dc(X)$), and they considered spaces with our property
(2), which they called DC-spaces.  With somewhat different terminology,
they prove what amounts to the fact that (6) implies (2), and that
(1) and (2) are equivalent when $X$ is a LOTS.

\bigskip

{\bf \S4. On Eberlein Compacta.}  Here we consider the properties
(1) -- (6) of \S3 in the case that $X$ is an Eberlein compact.
We already know by Theorem 3.3 that the stronger conditions
(4) or  (5) can hold only in the trivial case that the isolated
points of $X$ are dense in $X$; it is easy to see that (6)
holds iff $X$ is finite.  Thus, only (1), (2), and (3) are of
interest, and for these, the Eberlein compacta can be tailored to 
satisfy whatever we want.
The one we constructed in \S2 satisfies
(1), but not (2).  We now describe two modifications of this
construction, producing Eberlein compacta which satisfy
(3) but not (4) (this is easy), and then
(2) but not (3) (this requires more work).

For the first example, since we already know that (4) cannot hold
unless the isolated points are dense, it suffices to prove the following.

\medskip

{\bf 4.1. Theorem.}  There is a first countable uniform
Eberlein compact space $X$ such that $X$ has no isolated points and
$E_0(X,M) = C(X,M)$ for all metric spaces $M$.

{\bf Proof.}  Follow exactly the notation in $\S2$, so that $X$ will
be the $L_\omega$ there.  Choose sets $K_t$ for $t \in \cclw$ so that
conditions (Ra)(Rb)(Rc) of \S2.6 hold, so that all the lemmas through
Lemma 2.10 still apply.  But, replace (Rd)(Re) by
\smallskip

\vbox{  
\itemitem{\bf (Rf)}  Each $K_s$ is a singleton, and the 
$K_{s\alpha}$, for $\alpha\in\cc$, enumerate all the singletons
in $\hat K_s \backslash K_s$.
\smallskip
\noindent
As before, $U_t = X \cap (\pi_n^{-1}(K_t)\backslash K_t)$ where
$n = lh(t)$.
}

Now, fix $f \in C(L_\omega, M)$, where $M$ is metric.

{\bf Note:} $f$ is constant on $U_t$ for all but countably
many $t$.  If not, we could find an $s$ and an uncountable
$A \subseteq \cc$ such that $f$ is not constant on
$U_{s\alpha}$ for all $\alpha\in A$.
For $\alpha\in A$, let $K_{s\alpha} = \{\vec x_\alpha\}$, and choose
$\vec y_\alpha\in U_{s\alpha}$ such
that $f(\vec y_\alpha) \ne f(\vec x_\alpha)$.
Since the range of $f$
is compact, and hence second countable, we may,
as in the last paragraph of the proof of Lemma 2.4,
fix {\it distinct\/}
$p,q \in M$  and distinct $\alpha_n \in A$ ($n \in \omega$) such that
the $f(\vec x_{\alpha_n})$ converge to $p$ and
the $f(\vec y_{\alpha_n})$ converge to $q$.
Now, the points $\vec x_{\alpha_n}$
are in $K_s$,
which is compact metric, so, 
by passing to a subsequence, we may assume that the
$\vec x_{\alpha_n}$ converge to some point $\vec x \in K_s$.
Hence, in the weak topology, 
since  $\pi_{n+1}(\vec y_{\alpha}) = \vec x_\alpha$,
the $\vec y_{\alpha_n}$ converge to $\vec x$ also.
Applying $f$ to these sequences,
$f(\vec x) = p \ne q = f(\vec x)$, a contradiction.

It follows that $\Omega_f$ is dense in $X$, since every
nonempty open set in $X$ contains uncountably many $U_t$
(to see this, apply the above ``note'' and the fact that 
the co-zero sets of continuous functions form a basis for $X$).  \eop

\medskip

We remark that in the above ``note'', we used the same method to prove
$E_0(X)$ big as we used in Lemma 2.4
to prove $E_0(X)$ small; we have simply reversed the roles
of $f$ and $\pi$.

Also, it is possible to make the space of Theorem 4.1 zero-dimensional
by restricting the coordinates to lie in a Cantor set.  This would
not be possible for the spaces of \S2, or the space used
for Theorem 4.2(b) below.

Observe that in the proof of Theorem 4.1, 
the Hilbert space $\BB$ can be either complex or real,
since unlike in \S2, we no longer need the $\vec e_s$ direction 
to be two dimensional.  This holds in the next
construction as well, although we shall need that the base level
$L_0$ be infinite dimensional.

Also observe that if the
$K_{s\alpha}$ were not singletons, the above proof would 
establish a modified ``note'': for all but countably
many $t$, $f(\vec y) = f(\pi_n(\vec y))$ for all $\vec y \in U_t$.
This is the key to building a space satisfying (2) but not (3).
We shall make sure that $\hat K_s$ has ``large dimension'', so
that any real-valued function will be constant on many
subsets of $\hat K_s$, and these subsets will be the $K_{s\alpha}$;
this will ensure that $E_0(X,\RR) = C(X,\RR)$.  However,
if $M$ itself has ``large dimension'', then this argument will
fail, so that $E_0(X,M) \ne C(X,M)$.

The following definition and theorem pin down precisely for which
$M$ we can conclude $E_0(X,M)  = C(X,M)$ from 
$E_0(X,\RR) = C(X,\RR)$.  It suffices to consider only
compact $M$, since the range of each continuous map is compact.
Let $\FF_0$ be the collection of all zero or one point spaces.
For an ordinal $\alpha > 0$, let $\FF_\alpha$ be the class of all 
compact metric spaces $M$ such that there is a $\varphi \in C(M, [0,1])$,
with $\varphi^{-1}\{r\} \in \bigcup_{\delta < \alpha} \FF_\delta$ for every
$r \in [0,1]$.  So, for example, induction on $n \in \omega$ shows that
$[0,1]^n \in \FF_n$.  Then, if $M$ is the one-point compactification 
of the disjoint union of the $[0,1]^n$, we may let $\varphi$ map
$M$ to a simple sequence to conclude that $M \in \FF_\omega$.
Define $\FF = \bigcup_{\delta \in ON} \FF_\delta$, where $ON$
is the class of ordinals.  Actually, since 
every compact metric space has at most $\cc$ closed subspaces,
$\FF = \bigcup_{\delta < \cc^+} \FF_\delta$.

\medskip

\vbox{   
{\bf 4.2. Theorem.}
\item{a.} If $X$ is compact Hausdorff, $E_0(X,\RR) = C(X,\RR)$,
and $M \in \FF$, then $E_0(X,M) = C(X,M)$.
\item{b.}  There is a first countable uniform Eberlein compact $X$
such that $E_0(X,\RR) = C(X,\RR)$, but
for all compact metric spaces $M \notin \FF$,
$E_0(X,M) \ne C(X,M)$. }

\medskip 

\noindent
So, if we fix any compact metric space $M \notin \FF$,
we get an $X$ satisfying condition $(2)$ of
$\S3$, but not $(3)$.  Of course, we need to know that such
an $M$ exists, but that follows by a theorem of Levshenko.
There is a class of 
{\it strongly infinite dimensional\/} spaces which includes
the Hilbert cube, $[0,1]^\omega$.
Levshenko showed that if $M$ is a strongly infinite dimensional
compact metric space and $\varphi \in C(M,[0,1])$, then some
$\varphi^{-1}\{r\}$ is strongly infinite dimensional
(see [\AP]).  This gives us the following lemma.

\medskip

{\bf 4.3. Lemma.}
If $M \in \FF$, then $M$ is not strongly infinite dimensional.

{\bf Proof.} By induction on ordinals $\alpha$,
prove that every $M \in \FF_\alpha$ is not strongly infinite dimensional. \eop

\medskip

The definition of $\FF$ gives us the following 
easy inductive proof of Theorem 4.2(a).

\medskip

{\bf Proof of Theorem 4.2(a).}
Suppose that $M \in \FF_\alpha$, and suppose (inductively)
that the result holds for all $M' \in  \bigcup_{\delta < \alpha} \FF_\delta$.
Suppose $X$ is compact Hausdorff and $E_0(X,\RR) = C(X,\RR)$.
Fix $f \in C(X,M)$.  To prove $f \in E_0(X,M)$, we fix a nonempty
open $U \subseteq X$, and we produce a nonempty open
$V \subseteq U$ such that $f$ is constant on $V$.
Applying the definition of $\FF_\alpha$, fix
$\varphi \in C(M, [0,1])$ such that
for each $r \in [0,1]$ \ 
$\varphi^{-1}\{r\} \in \bigcup_{\delta < \alpha} \FF_\delta$.
Then $\varphi \circ f \in C(X,\RR) = E_0(X,\RR)$, so 
fix a nonempty  open set $W \subseteq U$ such that
$\varphi \circ f$ has some constant value $r$ on $\cl W$.
Now $\varphi^{-1}\{r\} \in \FF_\delta$, for some $\delta < \alpha$,
and  $E_0(\cl W,\RR) = C(\cl W,\RR)$ (by Lemma 3.1).
Applying the induction hypothesis, 
$f\rest \cl W \in C(\cl W, \varphi^{-1}\{r\} ) = E_0(\cl W,\varphi^{-1}\{r\})$,
so we may choose choose a nonempty open subset $V\subseteq W$ such
that $f\rest V$ is constant.   \eop

\medskip

To prove Theorem 4.2(b), we first prove some more lemmas about $\FF$.
Then, rather than construct a space $X$ which works for every compact
metric space $M \notin \FF$, we present Lemma 4.8, which allows us to
construct a separate $X_M$ for each $M$.  To construct each $X_M$,
we proceed as in \S2; that is, each $X_M$ will be an $L_\omega$,
constructed using somewhat modified conditions on the sets $K_t$.
We then glue these $X_M$ together to complete the proof of Theorem 4.2(b).

We begin with the lemmas about $\FF$.
First, another simple induction yields closure under subsets:

\medskip

{\bf 4.4. Lemma.}  If $M \in \FF$ and $H$ is a closed subset of $M$,
then $H \in \FF$.

\medskip

We also get closure under finite unions:

\medskip

{\bf 4.5. Lemma.}  Suppose that $M$ is compact metric and $M = H \cup K$,
where $H,K$ are closed subsets of $M$ and $H,K \in \FF$. 
Then  $M \in \FF$.

{\bf Proof.} Since $H$ is a closed $G_\delta$, fix
$\varphi \in C(M,  [0,1])$ such that $\varphi^{-1}\{0\} = H$.
Then, $\varphi^{-1}\{0\} \in \FF$.  For 
$r > 0$, we have $\varphi^{-1}\{r\} \subseteq K$,
so $\varphi^{-1}\{r\} \in \FF$ by Lemma 4.4.
So, $\varphi^{-1}\{r\} \in \FF$ for each  $r \in [0,1]$, which implies
that $M \in \FF$.  \eop

\medskip

Call $M$ {\it nowhere in $\FF$} iff $M$ is nonempty and
for each nonempty open $V \subseteq M$, we have
$\cl V \notin \FF$.
Note that such an $M$ has no isolated points, since
$\FF$ contains all one point spaces.

\medskip

{\bf 4.6. Lemma.}
If $M$ is a compact metric space and  $M \notin \FF$, 
then there is a closed set $K \subseteq M$ 
such that $K$ is nowhere in $\FF$.

{\bf Proof.}
Let $\UU = \{ U \subseteq M : U$ is open and
$\cl U \in \FF \}$, and let
$K = M \setminus \bigcup \UU$.

First, note that $K$ is nonempty:  If $K$ were empty, then, 
by compactness, $M$ would be covered by a finite subfamily
of $\UU$, which would imply $M \in \FF$ by Lemma 4.5.

To prove that $K$ is nowhere in $\FF$, it suffices
(by Lemma 4.4) to prove that 
$\cl {B(p,\epsilon)} \cap K \notin \FF$ whenever $p \in K$
and $\epsilon > 0$.  Note that $B(p,\epsilon)$ and its closure
are computed in $M$, not $K$.  Let $N = \cl {B(p,\epsilon)}$.
Fix $\varphi \in C(N,[0,1])$ such that
$\varphi^{-1}\{0\} = N \cap K$.
Since $B(p,\epsilon) \cap K$
is nonempty, $N \notin \FF$, so there must be some $r\in [0,1]$
such that $\varphi^{-1}\{r\} \notin \FF$.  
However, for $r > 0$,
$\varphi^{-1}\{r\}$ is compact and 
disjoint from $K$, so it is covered by a finite subfamily
of $\UU$,  and hence, as above, is in $\FF$.
So, $r$ must be $0$, so $N \cap K \notin \FF$.  \eop

\medskip

Let $\NNN(K)$ be the family of all compact $H \subseteq K$
such that $H$ is nowhere in $\FF$.  The following lemma is
trivial, given the above results, but we state it to
emphasize the abstract properties of our construction.

\medskip

{\bf 4.7. Lemma.} If $K$ is compact metric and
nowhere in $\FF$, then
\item{1.} $\NNN(K)$ is a family of nonempty closed subspaces of $K$.
\item{2.} $K \in \NNN(K)$.
\item{3.} For each $H \in \NNN(K)$ and each nonempty relatively open
$U \subseteq H$, there is an $L \in \NNN(K)$ with $L \subseteq U$.

\medskip

Most of the proof of Theorem 4.2(b) proceeds using just
the conclusion to Lemma 4.7, without any reference to $\FF$.
Note that if $K$ is a singleton, and $\NNN(K)$ is
redefined to be $\{K\}$, we also have the conclusion
to Lemma 4.7, and the proof of 4.2(b) then reproves Theorem 4.1.  

Now, as promised earlier, we present Lemma 4.8, which reduces
our construction to a modification of that of \S2.  

\medskip

{\bf 4.8. Lemma.}
If $\{X_{\alpha}: \alpha < \cc\}$ is a collection of nonempty
first countable uniform Eberlein compacta, then there is a
first countable uniform Eberlein compact space $X$, with
disjoint clopen subsets $J_{\alpha}$ homeomorphic to $X_{\alpha}$, 
such that $\bigcup_{\alpha \in \cc} J_{\alpha}$ is dense in $X$.

{\bf Proof.}
We may assume that each $X_{\alpha}$ is a weakly compact
subset of the closed unit ball of the Hilbert space $\BB_0$, and
that $\BB_0$ is a closed linear subspace of the Hilbert space $\BB$,
which contains unit vectors $\vec m_{\alpha} (\alpha \in \cc)$ and $\vec b$, 
all orthogonal
to each other and to $\BB_0$.

Let $r_{\alpha}$, for $\alpha \in \cc$, enumerate $(0,1)$.
Let $J_{\alpha} = X_{\alpha} + r_{\alpha}\vec b + \vec m_{\alpha}$.
Then $J_{\alpha}$ is homeomorphic to $X_{\alpha}$ (via translation).
Let $X$ be the union of the $J_{\alpha}$, 
together with all $r \vec b$ for $r \in [0,1]$.
Then $X$ is norm bounded (by $\sqrt3$), and is weakly closed, since any limit
of points in distinct $J_{\alpha}$ must be of the form $r\vec b$; the existence
of these limits also shows that the union of the $J_{\alpha}$ is dense in $X$.
The space $X$ is first countable by Lemma 2.3.  
To see that the $J_{\alpha}$ are disjoint and (weakly) clopen in $X$,
project along the $\vec m_{\alpha}$ direction. \eop

\medskip

We remark that translating along the $\vec b$ direction made 
$X$ first countable.  If we just let
$J_{\alpha} = X_{\alpha} + \vec m_{\alpha}$, and let $X$ be the union of the
$J_{\alpha}$ plus $\{\vec 0\}$,
then $X$ would be simply the one-point compactification of
the disjoint union of the $X_{\alpha}$.  Of course, we could build the
one-point compactification even if there are more
than $\cc \  X_{\alpha}$, in which
case one cannot make $X$ first
countable (by Arkhangel'ski\u\i's Theorem).

\medskip

We now construct the space  $X_M$.
Applying Lemma 4.6, let $K$ be a closed subset of $M$ which
is nowhere in $\FF$.
Let $\BB$ be a {\it real}
Hilbert space with an orthonormal basis consisting of unit
vectors $\{\vec e _ s : s \in \cclw \}
\bigcup \{\vec b _ i : i \in \omega \}$.
Let $\BB_n$ be the closed linear span of
$\{\vec e_s: lh(s) < n\} \cup \{\vec b _ i : i \in \omega \}$.
Since $\BB_0$ is infinite dimensional, we can embed $K$
in the first level of our space. To do so we replace condition (Ra)
of \S2.6 by the following:

\itemitem{\bf (Ra$\bf'$)} $K_{()}$ is a weakly compact subset of the closed
unit ball of $\BB_0$, and $K_{()}$ is homeomorphic to $K$.

\noindent
Actually, we could also make $K_{()}$ norm compact, but this is unnecessary.

Let $\pi_n$ be the perpendicular projection from $\BB$ onto $\BB_n$.
If $lh(s) = n$, let $D_s$ be the set of vectors of the
form $\vec v + \sum_{i \le n} r_i \vec e _ {s \rest i}$, 
where $\vec v \in K_{()}$ and each $|r_i| \le 2^{-i}$.
In particular, $D_{()}$ is homeomorphic to $K \times [-1,1]$.
As in \S2, the product with $[-1,1]$ allows us to make 
the $K_{\alpha}$ disjoint subsets of $D_{()}$.
As before, if $i \le n$, then $\pi_{i+1}(D_s) = D_{s \rest i}$.

We will choose the $K_t$ for $t \in \cclw$ so that
they satisfy condition {\bf (Ra$\bf'$)},
along with {\bf (Rb)} and  {\bf (Rc)} of \S2.6.  Now,  define 
$X_M = L_\omega$
to be the set of $\vec x \in \BB$ satisfying conditions (1), (2),
and (3) of \S2.7, along with
condition (0):  $\pi_0(\vec x) \in K_{()}$.

As before,  for $t \in \cclw$ and $n = lh(t)$, 
$U_t = L_\omega \cap (\pi_n^{-1}(K_t)\backslash K_t)$.
So $U_{()} =  L_\omega \setminus K_{()} = 
\{\vec x \in L_\omega : \vec x \cdot \vec e \ne 0\}$.
In this construction, we still have the levels $L_n = \pi_n(L_\omega)$,
with $L_0 = K_{()}$ and $L_1 = D_{()}$.
Now, elements of level $L_3 \setminus L_2$
are of the form
$\vec v + r_0 \vec e + r_1 \vec e_{\alpha} + r_2 \vec e_{\alpha\beta}$, 
where $0 < |r_i| \le 2^{-i}$ for each $i$,
$\vec v \in K_{()}$, 
$\vec v + r_0 \vec e \in K_\alpha$, and
$\vec v + r_0 \vec e + r_1 \vec e_\alpha \in K_{\alpha\beta}$.

\medskip

This $L_\omega$ still satisfies Lemmas $2.8$, $2.9$, and $2.10$, provided
we replace the bound in $2.9 (ii)$ by  ${7 \over 3}$.
The proofs are the same, except for the proof of $2.9 (iv)$,
where we join $\{\vec b _ i : i \in \omega \}$ to each  $C_{\vec x}$.

Now, we utilize $\NNN(K_{()})$ to  
choose the $K_s$.
Let $\vec u_{s\alpha}$ be of the
form $\vec u_s + r_{s\alpha} \vec e_s$,
with $\vec u_{()} = \vec 0$.
Choose the $K_s$ so that they satisfy,
in addition to {\bf (Ra$\bf'$)},
{\bf (Rb)}, and  {\bf (Rc)}, three more conditions:

\itemitem{\bf(Rg)}  Each $K_s$ is of the form $H_s + \vec u_s$,
where $H_s \in \NNN(K_{()})$.

\itemitem{\bf(Rh)}  For each $s$ and each $L \in  \NNN(K_{()})$ such that
$L \subseteq H_s$, $K_{s\alpha} = L + \vec u_{s\alpha}$ for
some $\alpha \in \cc$.

\itemitem{\bf(Ri)}  For each $s$ and each  nonempty relatively
open $V \subseteq \hat K_s$, there are uncountably many
$\alpha$ such that $K_{s\alpha} \subseteq V$.

\noindent
So, (Rg) says that each $K_s$ is a translate of a subset of $K_{()}$.
The $\hat K_s$ is defined precisely as in \S2, so that 
conditions  (Ra$'$),  (Rb),  (Rc) already imply that
$K_{s\alpha} \subseteq \hat K_s$.  Condition (Rg) guarantees that,
unlike in \S2, the projection $\pi_0 : K_t \to K_{()}$ is 1-1
for each $t$ (and its inverse is translation by $\vec u_t$). 
Using Lemma 4.7, it is easy to see
that conditions (Ra$'$),  (Rb),  (Rc), (Rg), (Rh), (Ri) can all be met.

If $f$ is a function on $L_\omega$ and $n = lh(t)$, we shall say
that $f$  is
$t$-{\it extension-constant\/} iff for all $\vec x \in K_t$
and all $\vec y \in L_\omega \cap \pi_n^{-1}\{\vec x\}$,
$f(\vec y) = f(\vec x)$.  By repeating the proof of the
``note'' in the proof of Theorem 4.1, we see the following:

\medskip

{\bf 4.9. Lemma.} If $M$ is metric and
$f \in C(L_\omega, M)$, then
$f$ is $t$-extension-constant
for all but countably many $t$.

\medskip
In the next lemma, we use condition (Ri) 
to show that the $U_t$ form a pi-base.

\medskip

{\bf 4.10. Lemma.} If $V$ is open and nonempty in  $L_\omega$,
then for some $t$, \  $U_t \subseteq V$.

{\bf Proof.} We may assume that
$V = \{\vec x \in L_\omega : f(\vec x) \ne 0\}$, where
$f \in C(L_\omega, \RR)$.   First fix $s$ such that 
$V \cap \hat K_s$ is nonempty, and then apply condition (Ri)
plus Lemma 4.9 to set $t = s\alpha$, where $\alpha$ is
chosen so that $K_{s\alpha} \subseteq V \cap \hat K_s$
and $f$ is $s\alpha$-extension-constant.  \eop

\medskip

In the case of Theorem 4.1, all the $K_t$ were singletons,
so ``$t$-extension-constant'' meant ``constant'', and the instance of Lemma 4.10
used there was simple enough that we omitted the proof of it.
In general, we cannot improve Lemma 4.9 to 
conclude that $f$ is constant on any open set.  
For example, the projection $\pi_0$ is 1-1 on each $K_t$, so
cannot be constant on
$K_t$ unless $K_t$ is a singleton. Applying Lemma 4.10, we get our
last lemma.

\medskip

{\bf 4.11. Lemma.} If 
$\NNN(K_{()})$ contains no singletons, then
$\pi_0 \in C(L_\omega, K_{()})$ and $\Omega_{\pi_0} = \emptyset$.

\medskip

Note, by condition (Ri), however, that
$\NNN(K_{()})$ contains no singletons iff no set in
$\NNN(K_{()})$ has any isolated points.  Of course, this is certainly 
true with $\NNN$ meaning ``nowhere in $\FF\,$''.  The specific features
of this $\NNN$ appear in the conclusion of our proof.

\medskip

{\bf Proof of Theorem 4.2(b).}  By Lemma 4.8 and the fact that
there are only $\cc$ compact metric spaces (up to homeomorphism),
it suffices to fix an $M \notin \FF$ and verify
that for the space $L_\omega$ constructed above,
$E_0(L_\omega,\RR) = C(L_\omega,\RR)$, but
$E_0(L_\omega,M) \ne C(L_\omega,M)$. 
Here, $L_\omega$ was constructed with $K_{()}$ homeomorphic
to a subset of $M$ which was nowhere in $\FF$, so that 
$E_0(L_\omega,M) \ne C(L_\omega,M)$ follows from Lemma 4.11.

Now, fix $f \in C(L_\omega,\RR)$.  In view of Lemma 4.10, to
prove that $f \in E_0(L_\omega,\RR)$, it suffices to
fix an $s$ and find a nonempty open $V \subseteq U_s$ on
which $f$ is constant.  By Lemma 4.10, fix $\alpha$ such
that $f$ is $s\alpha$-extension-constant.  By condition (Rg),
$K_{s\alpha} = H_{s\alpha} + \vec u_{s\alpha}$, where 
$H_{s\alpha} \in \NNN(K_{()})$.  Now, applying the properties
of $\NNN$, we can choose an $L \in  \NNN(K_{()})$ such that
$L \subseteq H_{s\alpha}$ and $f$ is constant on 
$L + \vec u_{s\alpha}$.  Applying condition (Rh)
to $s\alpha$, we can choose a $\beta$ such that
$H_{s\alpha\beta} = L + \vec u_{s\alpha\beta}$.
So let $V = U_{s\alpha\beta}$.  \eop

\bigskip

{\bf \S5. On Banach Spaces.}  In this section, we make a few remarks
on $E_0(X,M)$ in the case that $X$ is an arbitrary 
compact Hausdorff space and $M$ is a Banach space.
For definiteness, we take the scalar field to be $\RR$,
but all the results are unchanged if we replace $\RR$ by $\CC$.

First, as we 
have seen in \S3, there are many $X$ for which $E_0(X,M) = C(X,M)$.
For a given $X$, this can depend on $M$, but in view of \S4 and
the fact that every infinite dimensional Banach space contains
a homeomorphic copy of the Hilbert cube, there are only three possibilities:
\item{1.} $E_0(X,M) = C(X,M)$ for all Banach spaces $M$.
\item{2.} $E_0(X,M) = C(X,M)$ for all finite dimensional $M$, but not for any
infinite dimensional $M$.
\item{3.} $E_0(X,M) \ne C(X,M)$ for all Banach spaces $M$.

\noindent
Furthermore, there are Eberlein compact $X$ with no isolated points
realizing each
of these possibilities ((3) is trivial; see \S4 for (1) and (2)).

Second, in studying the properties of $E_0(X,M)$ 
as a normed linear space, we can isolate the two properties which
are of fundamental importance.  If $f,g_1,g_2 \in C(X,M)$, let us say that 
$f$ is {\it refined by\/} $g_1,g_2$  iff 
for all $x,y \in X$, if $g_1(x) = g_1(y)$ and $g_2(x) = g_2(y)$
then $f(x) = f(y)$.
A linear subspace $E \subseteq C(X,M)$ has the {\it refinement property\/}
iff for all $f,g_1,g_2 \in C(X,M)$, if
$g_1,g_2 \in E$ and $f$ is refined by $g_1,g_2$, then $f\in E$.
We say that $E$ has the {\it disjoint summation property\/} iff
whenever $\sum_{i\in\omega} f_i = f$ in $C(X,M)$, each $f_i \in E$,
and the sets $\{x: f_i(x) \ne 0\}$, for $i \in\omega$, are all disjoint,
then $f \in E$. 
The set of polynomial functions in $C([0,1],\RR)$ has
the disjoint summation property (trivially) but not the refinement property,
while the set of functions which are constant in some neighborhood
of ${1 \over 2}$ has the refinement
property but not the disjoint summation property.
Let us call $E$ a {\it nice\/} subspace of $C(X,M)$ iff $E$
has both properties.
Examples of nice $E$ are $E_0(X,M)$, $C(X,M)$, and the
space of all constant functions.  Or, one may fix any open $U \subseteq X$;
then $\{f \in C(X,M) : U \subseteq \overline\Omega_f\}$ is nice.
Another example is the functions of essentially countable range;
that is, let $\mu$ be a Baire measure on $X$, and
then let $D(X,M,\mu)$ be the set of $f\in C(X,M)$ such that
for some $\mu$-null-set $S \subseteq X$, $f(X\backslash S)$ is countable.
Another is the category analog of this --
the set of $f\in C(X,M)$ such that for some countable $P \subset M$, 
$\bigcup\{int(f^{-1}\{p\}) : p \in P\}$ is dense in $X$ ($int$ 
denotes ``interior'').

One advantage of studying nice $E$ is that we may restrict our
attention to the case where $E$ separates the points of $X$.
In general, given $E \subseteq C(X,M)$, we may
define an equivalence relation $\sim$ on $X$ by
$x \sim y$ iff $f(x) = f(y)$ for all $f\in E$.
Let $Y$ be the quotient, $X/\sim$; then
$Y$ is a compact Hausdorff space, and there is a canonical projection,
$\pi$, from $Y$ onto $X$.  Let
$E' = \{g \in C(Y,M) : g \circ\pi \in E\}$.  Then $E'$ is isometric
to $E$, and $E'$ separates the points of $Y$.  Further, both the
refinement property and the disjoint summation property 
are preserved here, so if $E$ is nice, then so is $E'$.

Some examples, when we start with $E = E_0(X,M)$:
For the spaces constructed in \S2:  If $X = L_\omega$, then $Y$
is a singleton.  If $X = L_2$, then $Y$ is obtained by collapsing
$L_1$ to a point.
In these two cases, $E' = E_0(Y,M)$, but this is not in general true.
For example, let $Q$ be any dense subset of $[0,1]$,
and form $X$ by attaching a copy of the $L_\omega$ of \S2 to each $q\in Q$,
where each copy goes off in some perpendicular direction.  There is then
a natural retraction, $r: X \to [0,1]$, and $E_0(X,\RR)$
consists of all functions of the form $f\circ r$, where $f\in C([0,1],\RR)$.
So, we may identify $Y$ with $[0,1]$ and $\pi$ with $r$, and
$E'$ is $C(Y,\RR)$, not $E_0(Y,\RR)$.

Third, we remark on some consequences of assuming that 
$E\subseteq C(X,M)$ has the refinement property.
If $\varphi\in C(M,M)$ and $f\in E$, then $\varphi\circ f \in E$
(since $\varphi\circ f$ is refined by $f,f$).
If $M = \RR$, and we view $C(X,M)$ as a Banach algebra
(under pointwise multiplication), then $E$ is a subalgebra.
More generally, if we fix any non-zero vector $\vec v\in M$, we may
let $\hat E \subseteq C(X,\RR)$ be the set of all 
$g \in C(X,\RR)$ such that the map $x \mapsto g(x) \vec v$ is in $E$.
Note that this does not depend on the $\vec v$ chosen, and
if $g \in \hat E$ and $f \in E$, then $g f \in E$.
Note also that $\hat E$ is nice.

It follows that 
if $E\subseteq C(X,M)$ has the refinement property and 
separates the points of $X$, then 
$E$ is dense in $C(X,M)$.  To see this, fix $f \in C(X,M)$.
If $M = \RR^n$, just apply the Stone-Weierstrass Theorem
to $f$ composed with the projections onto $n$ one-dimensional subspaces.
Then, for a general $M$, first approximate $f$ arbitrarily closely
by a map into a finite-dimensional subspace.

Actually, one can get more than just what is provided by
a simple application of the Stone-Weierstrass Theorem.
For example, we can arrange for the approximating function to
be identically zero wherever $f$ is zero:

\medskip

{\bf 5.1. Lemma.}  Suppose $E\subseteq C(X,M)$
has the refinement property
and separates the points of $X$.
Fix $f \in C(X,M)$ and fix  $\epsilon > 0$.  Then there is a
$g \in E$ with  $\|g - f\| \le \epsilon$,
$\|g\| = \|f\|$, and
$\|g(x)\| = \|f(x)\|$ for all $x$ such that $\|f(x)\|$ equals either
$0$ or $\|f\|$.

{\bf Proof.}  Assume $\epsilon < \|f\|$.
Fix $h \in E$ with  $\|h - f\| \le \epsilon /2$.
Then, let $\varphi : M \to M$ be any continuous map such that for all
$\vec v \in M$:  $\|\varphi(\vec v) - \vec v\| \le \epsilon /2$,
$\varphi(\vec v) = \vec0$ when $\|\vec v\| \le \epsilon / 2$,
and
$\varphi(\vec v) = \|f\|$ when $|\, \|\vec v\| - \|f\|\, | \le \epsilon / 2$
($\varphi$ can just move each $\vec v$ radially).  Then,
let $g = \varphi \circ h$. \eop

\medskip

Fourth, is $E_0(X,M)$ a Banach space?  Certainly it is in the
extreme cases where it is all of $C(X,M)$ and where
it contains only the constant functions.  To analyze the
general situation, we may, as pointed out above,
just consider the case where $E \subseteq C(X,M)$
is nice and separates the points of $X$.  Then, clearly, $E$ a Banach space in
the standard norm iff it is all of $C(X,M)$.  Furthermore, if
$E$	is not all of $C(X,M)$, then, following Bernard 
and Sidney [\Bera,\Sid], it is not even
{\it Banachizable\/}; that is, there is no norm  which
makes $E$ into a Banach space and gives $E$
a topology finer than the one inherited from $C(X,M)$.
In fact, every nice $E$ is {\it barreled\/}, which is
a stronger property.  There are a number of equivalents to being
barreled, discussed in [\Sid].  One is that for every linear space 
$\LL$ with $E \subseteq \LL\varsubsetneq  \overline E$,
$\LL$ is not Banachizable ($\overline E$ is the completion of
$E$; here, $\overline E =  C(X,M)$).
Another is the ``weak sequential property''
for $E$, which is the conclusion of the next Lemma; this is
a convenient way of establishing barreledness.
The proof of the next Lemma is very similar in spirit to
that of Theorem 2 of [\Sid], but we include it because
at first sight, the proof as stated in [\Sid] might appear to require
some additional assumptions about $E$ and $X$.
The two examples above of subspaces of $C([0,1],\RR)$ show that neither of the
two components of ``nice'', ``refinement property'' and 
``disjoint summation property'', is sufficient here.

\medskip

{\bf 5.2. Lemma.}  Let $X$ be compact and let $M$ be a Banach space.
Suppose that $E$ is a nice subspace of $C(X,M)$.
Let $\Lambda_n$, for $n \in \omega$, be in the dual space,
$E^*$.  Assume that for every $g \in E$,
$\Lambda_n(g) \to 0$.  Then $\sup_n\|\Lambda_n\| < \infty$.

{\bf Proof.}  As pointed out above, we may assume also
that $E$ is dense in $C(X,M)$,
so we may consider $\Lambda_n$ to be in $C(X,M)^*$.
Note that {\it if\/} $E = C(X,M)$,
the conclusion is immediate by the Banach-Steinhaus Theorem.
In any case, whenever $\HH$ is a closed
linear subspace of $C(X,M)$ such that
$\HH \subseteq E$,
$$
\sup\{|\Lambda_n(h)| : h \in \HH \cap \overline B(0, 1) \ \&\ 
n\in\omega \} < \infty\ \ . \eqno{(1)}
$$
Here, $\overline B(0, 1)$ is the closed unit ball of $C(X,M)$.
Now, assume that  $\sup_n\|\Lambda_n\| = \infty$.
We shall get a contradiction by applying (1).

For any $f \in C(X,M)$, let
$supt(f)$ be the closure of $\{x \in X : f(x) \ne \vec0\}$.
By compactness of $X$, we may fix a point $p$ such that
for all neighborhoods $V$ of $p$, 
$$
\sup\{|\Lambda_n(f)| : f \in \overline B(0, 1) \ \&\ n\in\omega \ \&\ 
supt(f) \subseteq V \} = \infty \ \ .
$$
By Lemma 5.1 (applied to $\hat E$ -- see above),
let $g \in \hat E$ be such that $\|g\| = 1$, $g(p) = 1$,
and $supt(g) \subseteq V$.
Then $\HH = \{g\vec v : \vec v \in M\}$ is a closed 
linear subspace of $C(X,M)$ (isometric to $M$) such that
$\HH \subseteq E$, so we may apply (1) above.
It follows, by considering functions of the form
$x \mapsto f(x) - g(x) f(p)$, that for all neighborhoods $V$ of $p$, 
$$
\sup\{|\Lambda_n(f)| : f \in \overline B(0, 1) \ \&\ n\in\omega \ \&\ 
supt(f) \subseteq V \ \&\ f(p)= \vec0  \} = \infty \ \ .
$$
Next, we show that for all neighborhoods $V$ of $p$, 
$$
\sup\{|\Lambda_n(g)| : g \in \overline B(0, 1)\cap E \ \&\ n\in\omega \ \&\ 
supt(g) \subseteq V\backslash\{p\} \} = \infty \ \ .
$$
To see this, fix $K > 0$, and then fix $n$ and 
$f \in C(X,M)$ such that $\|f\| \le 1$, $supt(f) \subseteq V$,
$f(p)= \vec0$, and $|\Lambda_n(f)| \ge 3K$.
Let $f' \in C(X,M)$ be such that $\|f'\| \le 1$, $supt(f') \subseteq V$,
$f'$ vanishes in some neighborhood of $p$, and
$\|f' - f\| \le K/\|\Lambda_n\|$.
Applying Lemma 5.1 to $f'$, let
$g \in E$ be such that $\|g\| \le 1$, $supt(g) \subseteq V$,
$g$ vanishes in some neighborhood of $p$, and
$\|g - f'\| \le K/\|\Lambda_n\|$.  Then 
$|\Lambda_n(g)| \ge K$.

Thus, we may inductively choose open neighborhoods $V_j$ of $p$,
$n_j \in\omega$, and $h_j \in E $  such that each
$\overline V_{j+1} \subseteq V_j$, 
$supt(h_j) \subseteq V_j \backslash \overline V_{j+1}$,
$\|h_j\| = 1$, and $|\Lambda_{n_j}(h_j)| \ge j$.
Let $\HH$ be the closed linear span in $C(X,M)$ of the $h_j$.
Since the $h_j$ are disjointly supported, 
$\HH \subseteq E$ (and $\HH$ is isometric to $c_0$),
so we have a contradiction to (1) above.  \eop

\medskip

Fifth, is $E_0(X,M)$ first category in itself?  We ask this because if
$E_0(X,M)$ is of second category, then Lemma 5.2 becomes trivial
by the Banach-Steinhaus Theorem.  Fortunately,
$E_0(X,M)$ is first category in  many cases; for example, when
$X$ contains a nonempty separable open subset with no isolated points
(see the proof of $(2)\imp(1)$ of Theorem 3.2).  In fact, as pointed
out by Bernard and Sidney, the original interest of $E_0(X)$ was
that it provided examples of first category normed linear spaces
which satisfy the Banach-Steinhaus Theorem, as well as a number
of other results usually proved by category arguments.
The following lemma describes some other situations in which
$E_0(X,M)$ is of first category.

\medskip

\vbox{
{\bf 5.3. Lemma.}  Let $X$ be compact and let $M$ be a Banach space.
Suppose that $E_0(X,M)$ is not a Banach space.  Then
$E_0(X,M)$ is of first category in itself
if either of the following hold:
\item{a.} $M$ is infinite dimensional.
\item{b.} $X$ is zero-dimensional.
}

{\bf Proof.}  First, as indicated above, we may pass to a
quotient and consider a nice $E \subseteq C(X,M)$ which is
dense in $C(X,M)$ but not all of $C(X,M)$;  of course, in (b),
this quotient operation is trivial.
Now, we need only show that $E$ is of first category in $C(X,M)$.

Whenever $H$ and $K$ are closed subsets of $X$,
let $U(H,K) = \{g \in C(X,M) : g(H) \cap g(K) = \emptyset\}$.
Note that $U(H,K)$ is always open in $C(X,M)$.  If $H$ and $K$
are disjoint, then either (a) or (b) guarantees that $U(H,K)$
is dense in $C(X,M)$.

Fix an $f \in C(X,M) \backslash E$.
Since $f(X)$ is second countable, there are closed $H_n, K_n \subseteq X$
for $n \in \omega$ such that each $H_n \cap K_n = \emptyset$,
and for all $x,y \in X$, if $f(x) \ne f(y)$, then for some $n$,
$x \in H_n$ and $y \in K_n$.  Let
$G = \bigcap_{n\in\omega}U(H_n,K_n)$.
Then $G$ is a dense $G_\delta$, and $f$ is refined by $g,g$ for all $g \in G$,
so $G$ is disjoint from $E$.  \eop

\medskip

The situation for finite dimensional $M$ seems more complicated.
We do not actually have an example of an 
$E_0(X,\RR)$ which is 
second category but not a Banach space, although it is easy to
produce a {\it consistent\/} example of this by forcing 
\hbox{[\Jech, \Kuna].}
In the ground model, $V$, let $X = L_\omega$ be the space
constructed in the proof
of Theorem 4.2(b), so $E_0(X,\RR) = C(X,\RR)$.
Let $V[G]$ add one Cohen real.  Then, in $V[G]$, 
$E_0(X,\RR)$ is of second category, since it contains the ground
model $C(X,\RR)$, which is of second category with this forcing.
However, in $V[G]$, $E_0(X,\RR)$ is not all of $C(X,\RR)$, since
$V[G]$ will contain a $g\in C(K_{()},\RR)$ which is 1-1 on
$K_{()} \cap V$; if $f = g \circ \pi_0 \in C(X,\RR)$, then 
$\Omega_f = \emptyset$.  To verify the details of this
construction, one must compare $X$ and $C(X,\RR)$ in
both models, $V$ and $V[G]$; this is described in \S3 of [\DK].

The following lemma yields a class of examples where $E_0(X,\RR)$ is of
first category.

\medskip

{\bf 5.4. Theorem.}  Let $M$ be any Banach space, and let
$X = \prod_{i\in\omega}X_i$, where each 
$X_i$ is compact Hausdorff and has more than one point.
Then $E_0(X,M)$ is of first category, and is dense in $C(X,M)$.

{\bf Proof.}  Let $P_n = \prod_{i= 0}^n X_i$,  and let
$\sigma_n$ be the projection from $X$ onto $P_n$.
Call a function $f$ on $X$ 
$n$-{\it supported\/} iff $f = g\circ \sigma_n$ for some
function $g$ on $P_n$.

To prove that $E_0(X,M)$ dense in $C(X,M)$, it is sufficient to
show that $E_0(X,\RR)$ separates points.  Fix two distinct points,
$x,y \in X$.  Since an infinite product has no isolated points,
we may assume (by partitioning the index set into infinitely many infinite
sets) that each $X_i$ has no isolated points.  We may also assume
that $\sigma_0(x) \ne \sigma_0(y)$.  We now produce an 
$f$ in $E_0(X, \RR)$ which separates $x,y$. 

Note that if $\sigma_n(\Omega_f) = P_n$ for all $n$,
then $\Omega_f$ will be dense.   To obtain this situation, 
we shall focus on the dyadic rationals.
Let $D_n = \{j\cdot 2^{-n} : 0 \le j \le 2^n\}$; so,
$D_0 = \{0,1\}$ and $D_1 = \{0, {1\over2}, 1\}$.
Inductively choose $f_n  \in C(X, [0,1])$ so that:

\itemitem{1.} $x \in int(f_0^{-1}\{0\})$ and $y \in int(f_0^{-1}\{1\})$.

\itemitem{2.} $f_n$ is $n$-supported.

\itemitem{3.} $\| f_{n+1} - f_n \| \le 2^{-n}$.

\itemitem{4.} $f_n^{-1}\{q\}  \subseteq f_{n+1}^{-1}\{q\} $
whenever $q \in D_n$.

\itemitem{5.} $ \bigcup \{ \sigma_n(int( f_{n+1}^{-1}\{q\})) : 
q \in D_{n+1}\} = P_n$.

\noindent
Let $f = lim_n f_n$.  This limit exists by (3).
$\sigma_n(\Omega_f) = P_n$ for all $n$ by (4)(5).
$f$ separates $x,y$ by (1).  Condition (2) allows the inductive
construction of $f_{n+1}$.

Now, we prove that $E_0(X,M)$ is of first category in $C(X,M)$.
For each $n$, let $U_n$ be the set of all $f \in C(X,M)$ 
such that for all $z \in P_n$, 
$f$ is not constant on $\{x \in X : \sigma_n(x) = z\}$.
Then $U_n$ is dense and open in $C(X,M)$, and $\Omega_f = \emptyset$
whenever $f \in \bigcup_{n\in\omega} U_n$.  \eop

\medskip

We remark that the space $D(X,M,\mu)$ defined above is 
always dense in $C(X,M)$ (by modifying the proof of the
Urysohn Separation Theorem), and is always of first category, except in
the trivial case that $\mu$ is a countable sum of point masses,
where $D(X,M,\mu) = C(X,M)$.

\bigskip
\bigskip
\bigskip

\vbox{
\centerline{\bf References}

\parskip = 3pt plus1pt minus.5pt
\frenchspacing
\parindent = 10 true mm

\item{[\AP]}  A. V. Arkhangel'ski\u\i  \ and  V. V. Fedorchuk,
{\it General Topology I, Basic Concepts and Constructions,
Dimension Theory}, Springer-Verlag, 1990.

\item{[\BHMWZ]}   A. Bella, A. Hager, J. Martinez, S. Woodward, and H. Zhou,
Specker spaces and their absolutes, I,  Preprint.

\item{[\BMW]}   A. Bella, J. Martinez and S. Woodward,
Algebras and spaces of dense constancies, Preprint.

\item{[\BRW]}  Y. Benyamini, M. E. Rudin, and M. Wage, 
Continuous images of weakly compact subsets of Banach spaces, 
{\it Pacific J. Math.}~70 (1977) 309-324.

\item{[\Ber]}  A. Bernard, Une fonction non Lipschitzienne peut-elle op\'erer
sur un espace de Banach de fonctions
non trivial?,  {\it J. Functional Anal.}~122 (1994) 451-477.

\item{[\Bera]}  A. Bernard, A strong superdensity property
for some subspaces of $C(X)$, 
Pr\'epubli\-cation de l'Institut Fourier,
Laboratoire de Math\'ematiques, 1994.

\item{[\BerSid]}  A. Bernard and S. J. Sidney,
Banach like normed linear spaces, Preprint, 1994.

\item{[\DK]} M. D\v zamonja and K. Kunen,
Properties of the class of measure
separable spaces, {\it Fund. Math.}, to appear.

\item{[\Hal]} P. R. Halmos, {\it Lectures on Boolean Algebras}, Van
Nostrand Reinhold, 1963.

\item{[\Jech]} T. Jech, {\it Set Theory},
Academic Press, 1978.

\item{[\Kuna]} K. Kunen, {\it Set Theory},
North-Holland, 1980.

\item{[\MW]}   J. Martinez and  S. Woodward,
Specker spaces and their absolutes, II,  Preprint.

\item{[\RuRu]}  M. E. Rudin and W. Rudin, 
Continuous functions that are locally constant on dense sets,
Preprint, 1994.

\item{[\Sid]}   S. J. Sidney, Some very dense subspaces
of $C(X)$, Preprint, 1994.

\item{[\Sik]}   R. Sikorski, {\it Boolean Algebras}, Springer-Verlag, 1964.

}

\bye